\documentclass{article}
\usepackage{amssymb}
\usepackage[centertags]{amsmath}

\begin{document}

\title{Basic theory of a kind of linear pantograph equations}

\author{Cheng-shi Liu\\Department of Mathematics\\Northeast Petroleum University\\Daqing 163318, China
\\Email: chengshiliu-68@126.com}

 \maketitle

\begin{abstract}
 By introducing a kind of special
 functions namely exponent-like function, cosine-like function and sine-like function,  we obtain explicitly the basic structures of solutions of initial value problem at the original point for this kind of linear pantograph equations. In particular,  we get the complete results on the existence, uniqueness and non-uniqueness of the initial value problems at a general point for the kind of
 linear pantograph equations.


\end{abstract}

\section{Introduction}
Usual functional differential equations such as $y'(x)=y(x-x_0)$
also namely time-delay differential equation have been expansively
studied from theories and applications[1]. Here, we consider another kind of functional differential equations with the form
$y'(x)=y(\alpha x)$ where $0<\alpha<1$, which is also named as pantograph equation. Due to the following several reasons, it is rather difficult to study them.
Firstly, even for the simplest homogenous linear pantograph $y'(x)=y(\alpha x)$, its solution can not be presented in terms of elementary functions. Secondly, even for the simplest non-homogenous linear pantograph equation $y'(x)=\lambda y(\alpha x)+q(x)$, we cannot solve it by the routine methods such as the variation of constant and the Laplace transformation. Thirdly, for the initial value problems of the linear pantograph equation, for the different choices of initial points, the existence and uniqueness of the solutions have crucial differences. In general, at a general point, the solution of the initial value problem perhaps doesn't exist, or is not unique. These facts are against intuition. Finally,  we know that the energy is a conservative quantity for
the second order vibration equation $y''(t)=-y(t)$, but for the
equation $y''(t)=-\alpha y(\alpha t)$, the energy is not conserved
since it is not invariant under the time translation. In fact,
letting $z(t)=y(t+T)$, we have
 $z''(t)=-\alpha y(\alpha t+\alpha T)\neq -\alpha z(\alpha t)=-\alpha y(\alpha t+ T)$.
 This implies that the original point of time has a special meaning and can not be taken
 arbitrarily.  Therefore, although there are many analogies between usual ordinary (and functional) differential  equations  and the pantograph equations, the most basic structures of the solutions of the linear pantograph equations are still open.

 The pantograph equations (including the case of $\alpha>1$) have been studied for a long time.
In 1940, Mahler[2] had introduced such type functional differential
equations in number theory. In 1971, Fox et al[3] and Ockendon et
al[4] proposed this kind of equations as the models to study some industrial problems. Kato and Mcleod [5]  studied
the asymptotic properties of the solution of the equation $y'(x)=ay(\lambda
x)+by(x)$.  Carr and Dyson [6,7] studied the related problems
on complex domain. In 1972, Morris, Feldstein and Bowen published an important paper[8] in which they obtained some very important results on the subject. In particular, they proved that the solution of the equation $y'(x)=-y(\alpha x)$ with $y(0)=1$ has an infinity of positive zeroes, which means that at such zero points, the solutions of initial value problem will be not unique or do not exist. It is Prof.Iserles who called this kind of functional differential equations with the multiplication delay as
 the pantograph equations[9].  Up to now, there have been a large number of papers
to study this kind of equations in theory and applications(see, for example,[9-23,25-28]). Among these,  Iserles[9,11] studied the generalized and nonlinear pantograph equations and gave some deep results, and  Derfel and Iserles[12] dealt with the equation on the complex
 plane. Iserles and Liu[18] studied the
 integro-differential pantograph equation. Mallet-Paret and Nussbaum[19] discussed the
 analyticity and non-analyticity of solutions. As application,  Van Brunt and Wake[21]
studied a model in cell growth. Because there exist in general no
solutions in terms of elementary functions and known special functions
even for the simplest equation $y'(x)=y(\alpha x)$, Iserles and
Liu[22] used the generalized hypergeometric functions to solve some
integro-differential pantograph equations. Feldstein, Iserles and Levin[23] studied the embedding of this kind of equations into the infinite-dimensional ODE systems.

Recently, Atiyah and
Moore discussed a kind of  functional differential equations such as
$y'(x)=k(y(x-x_0)+y(x+x_0))$ in the study of some problems arising
in fundamental physics [24]. This kind of delay
differential equations  can be solved by elementary functions.
Atiyah and Moore pointed out that relativistic invariance implied
that one must consider both advanced and retarded terms in the
equations, and they named them as shifted equations and showed that
the shifted Dirac equation had some novel properties and a tentative
formulation of shifted Einstein-Maxwell equations naturally
incorporates a small but nonzero cosmological constant. Following
Atiyah and Moore[24], Kong and Zhang [25] also studied the
pantograph type equations.
If we notice that there is an original point of time
 for our universe which arises from the big bang, this kind of
 equations will become a possible mathematical tool to describe the
 corresponding physics. Another possible application of pantograph equations is elastic-plastic
 mechanics[29] in which the strain of the
 material depends on substraction delay type memory of the stress.
 Indeed, for instance, for an elastic-plastic spring,
 a memory effect can in general be represented in terms of a delay
 function $h(t)$ such that the force at the time $t$ depends on the displacement at the time
 $t-h(t)$. For simplicity, a routine way is to take a constant delay function
 $h(t)=h_0$ which leads to a usual  delay differential equation. But there exist some
 weaknesses for this choice. For example, the memory in $t<h_0$ can not be
 considered. Moreover, the fixed $h_0$ is only suitable for short
 period memory. For long time memory, a reasonable memory function
 should be a real function depending on a time variable $t$. The
 simplest choice is $h(t)=\beta t$ with $0<\beta<1$, and hence
 $t-h(t)=(1-\beta)t=\alpha t$ with $\alpha=1-\beta$. This treatment
 leads naturally to a multiplication delay differential equation. Therefore,
  if replacing the substraction delay by multiplication delay, the
 corresponding vibration equation will be $y''(t)=-Ky(\alpha t)$
 where $0<\alpha<1$. Some new applications can be found in [26-28].

 In the present paper,  we introducing a kind of special functions
  namely exponent-like function, cosine-like function and sine-like function, and explicitly solve
  the initial value problems of these equations including  the first order linear pantograph equation, the high order linear pantograph equation,
  the system of linear pantograph equations,  and the boundary
  value problem of the second order pantograph equation.
  We give some important and interesting properties of these  special functions.
  By these special functions, the structures of the solutions of these linear pantograph equations are recovered and obtained  clearly. In particular, we obtain the basic results on the existence, uniqueness and non-uniqueness of the solutions of initial value problem at a general point for linear pantograph equations. These results show some crucial differences between the usual linear differential equations and linear pantograph equations.

This paper is organized as follows. In section 2, as preliminary,
the existence and uniqueness theorems for the initial value problem at the original point of
  linear pantograph  equations are listed. Section 3 is a key section of the paper,
 in which three special functions
 namely exponential-like function $E_{\alpha}(x)$, sine-like function $S_{\alpha}(x)$ and cosine-like
 function $C_{\alpha}(x)$ are introduced and some properties are obtained by detailed analysis.
  These functions
 are the foundation of solving and studying the linear pantograph equations. In section 4, we give the basic theorems on linear pantograph equations. The non-homogenous pantograph equation
 $y'(t)=py(\alpha t)+q(x)$ is solved by the power series method and a useful formula is
 given. This is a key result by which the structure of the solutions
 of the system of linear pantograph equations is obtained.
 In section 5, the solutions of the second order pantograph equation $y''(x)+py'(\alpha x)+ qy(\alpha^2 x)=0$ are
 classified and its initial value problem at the original point is solved. In section 6, the system of the  linear pantograph equations is solved in details and the structure of the solutions is given. In section 7, an operator method is used to solve the high order linear pantograph equations and the corresponding formula of the solution is obtained. In section 8, we consider the initial value problem at a general point for linear pantograph equations, and obtain the complete results on the existence and uniqueness and non-uniqueness of solutions.
  In section 9, the boundary value problem of the second order linear pantograph equation is solved and the
  corresponding eigenvalues and eigenfunctions are given.

\section{Preliminary: existence and uniqueness theorems at original point}

Since the existence and uniqueness theorems are the foundations of
the further studies, we give them as the preliminary in this
section. For the following initial value problem of the first order pantograph equation
\begin{equation}
y'(x)= f(y(\alpha x)),
\end{equation}
\begin{equation}
y(0)=y_0,
\end{equation}
where $0<\alpha<1$, if $f$ satisfies the Lipschitz condition, the
existence and uniqueness of the local solution can be proven by the
Banach fixed point theorem. In the paper, we only consider the
linear equation, and hence we need not this result.

 By the standard method in usual ODE theory[30], we can easily prove the following theorems (see the theorem 1 in [5]).

\textbf{Theorem 2.1}. There is a unique analytic solution for the initial value problem
\begin{equation}
\frac{\mathrm{d}X(t)}{\mathrm{d}t}=AX(\alpha t)+b(t),
\end{equation}
\begin{equation}
X(0)=X_0,
\end{equation}
where $X(t)=(x_1(t),\cdots,x_n(t))^T$ is a $n$ dimensional column vector, $b(t)$ is a known column vector function and $A=(a_{ij})_{n\times n}$ is a constant matrix.

\textbf{Theorem 2.2}. If $X_1(t),\cdots,X_n(t)$ are $n$ linear independent basic solutions of the homogenous equation of Eq.(3), and $X^*(t)$ is a special solution of the non-homogenous equation, then the general solution of non-homogenous equation can be represented by linear combination of these basics solutions adding a special solution, i.e.,
\begin{equation}
X(t)=c_1X_1(t)+\cdots+c_nX_n(t)+X^*(t).
\end{equation}

For example, for the  first order linear pantograph equation
\begin{equation}
y'(x)= \beta y(\alpha x),
\end{equation}
\begin{equation}
y(0)=y_0,
\end{equation}
or the second equation
\begin{equation}
y''(x)= \beta y(\alpha x),
\end{equation}
\begin{equation}
y(0)=y_0, y'(0)=v_0,
\end{equation}
there is the unique analytic solution.

\textbf{Corollary 2.1}. If $y_1(x)$ and $y_2(x)$ are two linear
independent solutions of the equation $y''(x)=\beta y(\alpha x)$
with initial values $y_1(0)=1, y'(0)=0$ and $y_2(0)=0, y'(0)=1$
respectively, the general solution can be represented by
\begin{equation}
y(x)=c_1y_1(x)+c_2y_2(x),
\end{equation}
where $c_1$ and $c_2$ are toe arbitrary constants.

 \textbf{Proof}. Firstly, by initial conditions, it is easy to see that these two solutions $y_1$ and $y_2$ are linear independent.
 Let $y(x)$ be a solution  with initial conditions $y(0)=c_1$ and $y'(0)=c_2$, and take $z(x)=c_1y_1(x)+c_2y_2(x)$.
 Then $y(x)$ and $z(x)$ have the same initial values. From the unique theorem 2.1, it follows that $y(x)=z(x)$.

\textbf{Theorem 2.3}. For the $n$-th order linear pantograph equation
\begin{equation}
y^{(n)}(x)+p_{n-1}y^{(n-1)}(\alpha x)+p_{n-2}y^{(n-2)}(\alpha^2
x)+\cdots+p_1y'(\alpha^{n-1} x)+p_0y(\alpha^n x)=f(t),
\end{equation}
\begin{equation}
y(0)=c_1, y'(0)=c_2,\cdots,y^{(n-1)}(0)=c_{n-1},
\end{equation}
where $f(t)$ is a known analytic function,  there is a unique analytic solution, and the general solution of the Eq.(11) can be represented by linear combination of the basics solutions of homogenous equation adding a special solution of non-homogenous equation.

\textbf{ Proof}. letting $x_1(t)=y(t), x_2(t)=x'_1(\frac{t}{\alpha}),\cdots, x_n(t)=x'_{n-1}(\frac{t}{\alpha})$,  the equation (11) is transformed to the first order linear equations system.
By theorem 2.2,  the proof is completed.

We can easily prove the following theorem.

\textbf{Theorem 2.4}. For the $n$-th order linear pantograph equation (11),
if $y_k(x)$ is the special solution of the equation (11) with $f(t)=f_k(t)$, then
\begin{equation*}
y(x)=y_1(x)+\cdots+y_m(x)
\end{equation*}
is a special solution of Eq.(11)with $f(t)=\sum_{k=1}^m f_k(t)$.

 \section{Exponent-like function $E_{\alpha}(x)$, cosine-like function $C_{\alpha}(x)$ and sine-like function $S_{\alpha}(x)$}

 In the section, we introduce three special functions which are the
 nontrivial generalizations of the usual exponential
 function, cosine function and sine function, and have some interesting
 properties such as infinite addition formulae. This section plays an important role in the paper since these special functions are the key
 mathematical tools and the foundation of all following studies.

 We denote $E_{\alpha}(x)$ the unique analytic solution of the following initial value problem
 \begin{equation}
 y'(x)=y(\alpha x),
 \end{equation}
 \begin{equation}
 y(0)=1.
 \end{equation}
  Then its power series expansion is given as follows
 \begin{equation}
 E_{\alpha}(x)=\sum_{n=0}^{+\infty} \alpha^{\frac{n(n-1)}{2}}\frac{x^n}{n!}.
 \end{equation}
Replacing $x$ by $ix$ yields
 \begin{equation}
 E_{\alpha}(ix)= C_{\alpha}(x)+iS_{\alpha}(x),
 \end{equation}
where $i^2=-1$ and
\begin{equation}
 C_{\alpha}(x)=\sum_{n=0}^{+\infty}(-1)^n \alpha^{n(2n-1)}\frac{x^{2n}}{(2n)!},
 \end{equation}
 \begin{equation}
 S_{\alpha}(x)=\sum_{n=0}^{+\infty}(-1)^n \alpha^{n(2n+1)}\frac{x^{2n+1}}{(2n+1)!}.
 \end{equation}
 It is also easy to see that
 \begin{equation}
 C_{\alpha}(x)=\frac{E_{\alpha}(ix)+E_{\alpha}(-ix)}{2},
 \end{equation}
 and
 \begin{equation}
 S_{\alpha}(x)=\frac{E_{\alpha}(ix)-E_{\alpha}(-ix)}{2i}.
 \end{equation}

We call these three functions respectively the exponent-like function $E_{\alpha}(x)$, cosine-like function $C_{\alpha}(x)$ and sine-like function $S_{\alpha}(x)$. When $\alpha=1$, they all become the usual exponential function, cosine function and sine function. It is easy to show that these three special functions are analytic functions on whole complex plane, that is, they are entire functions, and $C_{\alpha}(x)$ is an even function and $S_{\alpha}(x)$ is an odd function, with $C_{\alpha}(0)=1$, $C'_{\alpha}(0)=0$, $S_{\alpha}(0)=0$ and $S'_{\alpha}(0)=1$. Further, we can easily prove that the orders of them are zero. In particular, we can easily prove the following important properties.

\textbf{Remark 3.1}. The power series form of the function $E_{\alpha}(x)$ had been obtained in some papers such as [5,9,16]. Here, I give it a suitable name and symbol so that it can be studied as an independent mathematical object and used as a convenient mathematical tool.

\textbf{Proposition 3.1}. For the first order derivative, we have
\begin{equation}
E'_{\alpha}(x)=E_{\alpha}( \alpha x),
\end{equation}
\begin{equation}
C'_{\alpha}( x)=- S_{\alpha}( \alpha x),
\end{equation}
\begin{equation}
S'_{\alpha}( x)= C_{\alpha}( \alpha x),
\end{equation}
where prime means the derivative with respect to $x$.

\textbf{Proposition 3.2}. For the second order derivative, we have
\begin{equation}
C''_{\alpha}(x)=- \alpha C_{\alpha}( \alpha^2 x),
\end{equation}
\begin{equation}
S''_{\alpha}( x)=-\alpha S_{\alpha}(\alpha^2 x),
\end{equation}
where prime means the derivative with respect to $x$.

Furthermore, we have the following addition formulae.

\textbf{Proposition 3.3} (Addition formulae). The exponent-like function $E_{\alpha}(x)$, cosine-like function $C_{\alpha}(x)$ and sine-like function $S_{\alpha}(x)$ satisfy the addition formulae
\begin{equation}
E_{\alpha}(x+y)= \sum_{n=0}^{+\infty} \alpha^{\frac{n(n-1)}{2}}\frac{x^n}{n!}E_{\alpha}(\alpha^ny)=\sum_{n=0}^{+\infty} \alpha^{\frac{n(n-1)}{2}}\frac{y^n}{n!}E_{\alpha}(\alpha^nx);
\end{equation}
\begin{equation}
C_{\alpha}(x+y)= \sum_{n=0}^{+\infty} \alpha^{n(2n-1)}\frac{(-1)^nx^{2n}}{(2n)!}C_{\alpha}(\alpha^{2n}y)-\sum_{n=0}^{+\infty} \alpha^{n(2n+1)}\frac{(-1)^nx^{2n+1}}{(2n+1)!}S_{\alpha}(\alpha^{2n+1}y);
\end{equation}
\begin{equation}
S_{\alpha}(x+y)= \sum_{n=0}^{+\infty} \alpha^{n(2n-1)}\frac{(-1)^nx^{2n}}{(2n)!}S_{\alpha}(\alpha^{2n}y)+\sum_{n=0}^{+\infty} \alpha^{n(2n+1)}\frac{(-1)^nx^{2n+1}}{(2n+1)!}C_{\alpha}(\alpha^{2n+1}y).
\end{equation}

\textbf{Remark 3.2}. When $\alpha=1$, these formulae are just the
usual addition formulae $\exp(x+y)=\exp(x)\exp(y),
\cos(x+y)=\cos(x)\cos(y)-\sin(x)\sin(y)$ and $
\sin(x+y)=\sin(x)\cos(y)+\cos(x)\sin(y)$.

\textbf{Proof}. Since the corresponding series are uniform convergent, the following summations can exchange orders.  By direct computation, we have

\begin{equation*}
E_{\alpha}(x+y)= \sum_{n=0}^{+\infty} \alpha^{\frac{n(n-1)}{2}}\frac{1}{n!}(x+y)^n
\end{equation*}
\begin{equation*}
= \sum_{n=0}^{+\infty} \alpha^{\frac{n(n-1)}{2}}\frac{1}{n!}\sum_{k=0}^n\frac{n!}{k!(n-k)!}x^ky^{n-k}
\end{equation*}
\begin{equation*}
= \sum_{n=0}^{+\infty} \alpha^{\frac{n(n-1)}{2}}\frac{x^n}{n!}\sum_{m=0}^{+\infty}\alpha^{\frac{m(m-1)}{2}}\frac{y^m}{m!}\alpha^{mn}
\end{equation*}
\begin{equation*}
= \sum_{n=0}^{+\infty} \alpha^{\frac{n(n-1)}{2}}\frac{x^n}{n!}E_{\alpha}(\alpha^ny)
\end{equation*}
\begin{equation*}
=\sum_{n=0}^{+\infty} \alpha^{\frac{n(n-1)}{2}}\frac{y^n}{n!}E_{\alpha}(\alpha^nx).
\end{equation*}
From the first formula (26), we have
\begin{equation}
C_{\alpha}(x+y)+C_{\alpha}(x-y)= 2\sum_{n=0}^{+\infty} \alpha^{n(2n-1)}\frac{(-1)^nx^{2n}}{(2n)!}C_{\alpha}(\alpha^{2n}y);
\end{equation}
\begin{equation}
S_{\alpha}(x+y)+S_{\alpha}(x-y)=2\sum_{n=0}^{+\infty} \alpha^{n(2n+1)}\frac{(-1)^nx^{2n+1}}{(2n+1)!}C_{\alpha}(\alpha^{2n+1}y);
\end{equation}
\begin{equation}
C_{\alpha}(x+y)-C_{\alpha}(x-y)= -2\sum_{n=0}^{+\infty} \alpha^{n(2n+1)}\frac{(-1)^nx^{2n+1}}{(2n+1)!}S_{\alpha}(\alpha^{2n+1}y);
\end{equation}
\begin{equation}
S_{\alpha}(x+y)-S_{\alpha}(x-y)= 2\sum_{n=0}^{+\infty} \alpha^{n(2n-1)}\frac{(-1)^nx^{2n}}{(2n)!}S_{\alpha}(\alpha^{2n}y).
\end{equation}
By the above these formulae, we obtain the last two formulae (27) and (28). The proof is completed.

\textbf{Remark 3.3}.  We define $y=L_{\alpha}(x)$ to be the inverse function of $x=E_{\alpha}(y)$ and call it the logarithm-like function. Therefore, we have
\begin{equation*}
\frac{\mathrm{d}y(x)}{\mathrm{d}x}=\frac{1}{x(\alpha y)}.
\end{equation*}
By letting $1+x=E_{\alpha}(y)$ and using $\frac{\mathrm{d}}{\mathrm{d}x}=\frac{\mathrm{d}y}{\mathrm{d}x}\frac{\mathrm{d}}
{\mathrm{d}y}=\frac{1}{E_{\alpha}(\alpha y)}\frac{\mathrm{d}}{\mathrm{d}y}$, we derive the power series expansion of $L_{\alpha}(1+x)$ as follows
\begin{equation}
L_{\alpha}(1+x)=x-\frac{\alpha }{2!}x^2+\frac{\alpha^2(3-\alpha) }{3!}x^3-\frac{\alpha^3(\alpha^3-6\alpha+11) }{4!}x^4+\cdots
\end{equation}
If $\alpha=1$, it is just the power series of the usual logarithm function $\ln(1+x)$.

On the zeros of cosine-like function $C_{\alpha}(x)$ and sine-like function $S_{\alpha}(x)$, we have the following theorems.

\textbf{Theorem 3.1}. For $0<\alpha\leq 1$, the cosine-like function $C_{\alpha}(x)$ and sine-like function  $S_{\alpha}(x)$ have respectively an infinity of real zeroes.

\textbf{Proof}. We only need to prove the result of $C_{\alpha}(x)$. Since  $C_{\alpha}(x)$ satisfies the equation
\begin{equation*}
y''(x)=-\alpha y(\alpha^2 x),
\end{equation*}
with $y(0)=1, y'(0)=0$, this means that the graph of $C_{\alpha}(x)$ starts at the point $(0,1)$ and moves to the right with slope beginning at zero. By the equation itself, we know that $y''(x)=-\alpha y(\alpha^2 x)$, so when the curve is above the $x$ axis, $C''_{\alpha}(x)$ is negative that increases as the curve decreases, and hence the curve  $C_{\alpha}(x)$ bends down and crosses the $x$ axis at some point $\eta_0$. Then at the point $ \frac{\eta_0}{\alpha^2}$ we have $C''_{\alpha}(\frac{\eta_0}{\alpha^2})=0$ so the curve has an inflection point, and then the curve goes down continuously to the local lowest point $\alpha\rho_1$ where $\rho_1$ is the first positive zero point of $S_{\alpha}(x)$ which exists by a similar discussion with $C_{\alpha}(x)$, and then
the curve bends up and crosses $x$ axis at some point $\eta_1$. This process will be continuous forever and gives an infinity of zeros. The proof is completed.

Next we give other properties of $C_{\alpha}(x)$ and  $S_{\alpha}(x)$. Denote $\rho_0=0, \pm\rho_1,\cdots,\pm\rho_n,\cdots$ and $ \pm\eta_1,\cdots,\pm\eta_n,\cdots$ as the zeros of  $S_{\alpha}(x)$ and $C_{\alpha}(x)$ respectively. It is easy to see that $S_{\alpha}(\rho_n)=0$ implies $C'_{\alpha}(\frac{\rho_n}{\alpha})=0$ and vice versa. $C_{\alpha}(x)$ is convex on the intervals $[\eta_{2k}/\alpha,\eta_{2k+1}/\alpha]$,  and $C_{\alpha}(x)$ is concave on the intervals $[\eta_{2k+1}/\alpha,\eta_{2k+2}/\alpha]$. The similar results hold for  $S_{\alpha}(x)$. Furthermore, we have the following theorem.

\textbf{Theorem 3.2}.  All real zeros of $C_{\alpha}(x)$ and $S_{\alpha}(x)$ are alternative each other, and their positive zeros satisfy
\begin{equation}
0=\rho_0<\alpha \eta_1<\eta_1<\alpha\rho_1<\rho_1<\alpha\eta_2<\eta_2<\alpha\rho_2<\rho_2<\cdots.
\end{equation}

\textbf{Proof}. Since $S_{\alpha}(\rho_n)=S_{\alpha}(\rho_{n+1})=0$, by Roll's theorem, there exists a point $r_n\in(\rho_n,\rho_{n+1})$, such that $S'_{\alpha}(r_n)=0$, that is
\begin{equation*}
C_{\alpha}(\alpha r_n)=S'_{\alpha}(r_n)=0.
\end{equation*}
It follows that $\alpha r_n$ is the $(n-1)$-th positive zero of $C_{\alpha}(x)$, that is, $\eta_{n-1}=\alpha r_n$, and then,
\begin{equation*}
\alpha\rho_{n-1}<\eta_n<\alpha\rho_{n}.
\end{equation*}
Similarly, we have
\begin{equation*}
\alpha\eta_n<\rho_n<\alpha\eta_{n+1}.
\end{equation*}
The proof is completed.

We can give another proof by using the theorem about the relation of entire function and its derivative. We omit it.

\textbf{Theorem 3.3} (The comparison theorem of the first positive zero). If $y(x)$ is the solution of equation
\begin{equation*}
y''(x)=-k y(\alpha^2 x),
\end{equation*}
with $k>0$ and initial conditions $y(0)=1, y'(0)=0$ or $y(0)=0, y'(0)=1$, then the first positive zero $x_0$ is a decreasing function as a function of $k$ with $\alpha$ fixed.

\textbf{Proof}. We only consider the case  $y(0)=1, y'(0)=0$. Let $y(x)$ and $z(x)$ be respectively solutions of equations
\begin{equation*}
y''(x)=-k_1 y(\alpha^2 x),
\end{equation*}
\begin{equation*}
z''(x)=-k_2 z(\alpha^2 x),
\end{equation*}
with $0<k_1<k_2$ and the same initial conditions. It is easy to see that $z(x)=y(\sqrt{\frac{k_2}{k_1}}x)$ is the solution of the second equation with the initial conditions $z(0)=1, z'(0)=0$. We assume that $x_0$ is the first positive zero of $y(x)$, then $x_1=\sqrt{\frac{k_1}{k_2}}x_0<x_0$ is the first positive zero of $z(x)$ since $z(x_1)=y(x_0)=0$. The conclusion in the case $y(0)=0, y'(0)=1$ can be proven by similar method. The proof is completed.

\textbf{Theorem 3.4}. At any ray, $S_{\alpha}(z)$ and  $C_{\alpha}(z)$ are unbounded. As a special case, $S_{\alpha}(x)$ and  $C_{\alpha}(x)$ are unbounded functions in real axis.

\textbf{Proof}. Since the orders of $S_{\alpha}(z)$ and  $C_{\alpha}(z)$ are zero, by Phragm\'{e}n-Lindel\^{o}f theorem see, [31]), we get the conclusion. The proof is completed.

\textbf{Theorem 3.5}. $S_{\alpha}(z)$ and  $C_{\alpha}(z)$ have only real zeros.

\textbf{Proof}. We only consider $C_{\alpha}(z)$. Firstly, we know that the function can be written as
\begin{equation*}
 C_{\alpha}(x)=\sum_{n=0}^{+\infty} \alpha^{2n^2}\frac{(-(\alpha^{-\frac{1}{2}}x)^2)^n}{(2n)!}.
 \end{equation*}
Denote $z=-(\alpha^{-\frac{1}{2}}x)^2$, then
\begin{equation*}
 C_{\alpha}(x)=G(z)=\sum_{n=0}^{+\infty} \alpha^{2n^2}\frac{z^n}{(2n)!}.
 \end{equation*}
According to the power series of cosine function, we know that
\begin{equation*}
 F(z)=\sum_{n=0}^{+\infty}\frac{z^n}{(2n)!},
 \end{equation*}
has only real zeros, and its infinite product is
 \begin{equation*}
 F(z)=\prod_{n=0}^{+\infty}(1+\frac{z}{z_n}),
 \end{equation*}
 where $z_n=(n+\frac{1}{2}\pi)^2$. In addition, $\{q_n=\alpha^{2n^2}\}_{n=0}^{+\infty}$ is a multiplier sequence(see, P\'{o}lya and Schur [32]). Therefore, by the Laguerre's theorem[33], $C_{\alpha}(z)$ only has real zeros which number is infinite. Similarly, we can get the same conclusion on $S_{\alpha}(z)$. The proof is completed.

\textbf{Proposition 3.4}. $\rho_m$ and $\eta_m$ satisfy the following identity,
\begin{equation*}
 \sum_{n=0}^{+\infty} \alpha^{n(2n-1)}\frac{(-1)^n(\rho_m-\eta_m)^{2n}}{(2n)!}S_{\alpha}(\alpha^{2n}\eta_m)
\end{equation*}
\begin{equation}
 +\sum_{n=0}^{+\infty} \alpha^{n(2n+1)}\frac{(-1)^n(\rho_m-\eta_m)^{2n+1}}{(2n+1)!}C_{\alpha}(\alpha^{2n+1}\eta_m)=0.
\end{equation}
\textbf{Proof.} Taking $x=\rho_m-\eta_m$ and $y=\eta_m$ in the formula (28) gives the result. The proof is completed.

This is a complicated relation between these zeros. Naturally, an interesting  problem is to study whether there exists a simple relation between the $n$-th zeroes $\rho_m$ and $\eta_m$ of $S_{\alpha}(x)$ and  $C_{\alpha}(x)$.
I leave it as an open problem.

On the integrals of $C_{\alpha}(x)$ and  $S_{\alpha}(x)$, we have the following results.

\textbf{Proposition 3.5}.
\begin{equation}
\int_{\frac{\rho_n}{\alpha}}^{\frac{\rho_{n+1}}{\alpha}}C_{\alpha}(\alpha^2 x)\mathrm{d}x=\int_{\alpha\rho_n}^{\alpha\rho_{n+1}}C_{\alpha}(x)\mathrm{d}x=0,
\end{equation}
\begin{equation}
\int_{\frac{\eta_n}{\alpha}}^{\frac{\eta_{n+1}}{\alpha}}S_{\alpha}(\alpha^2 x)\mathrm{d}x=\int_{\alpha\eta_n}^{\alpha\eta_{n+1}}S_{\alpha}(x)\mathrm{d}x=0.
\end{equation}

\textbf{Proof}. By $C''_{\alpha}(x)=-\alpha C_{\alpha}(\alpha^2 x)$ and $C'_{\alpha}(\frac{\rho_n}{\alpha})=0$, we have
\begin{equation*}
\int_{\frac{\rho_n}{\alpha}}^{\frac{\rho_{n+1}}{\alpha}}C_{\alpha}(\alpha^2 x)\mathrm{d}x=-\frac{1}{\alpha}\int_{\frac{\rho_n}{\alpha}}^{\frac{\rho_{n+1}}{\alpha}}C''_{\alpha}(x)\mathrm{d}x
=-\frac{1}{\alpha}(C'_{\alpha}(\frac{\rho_{n+1}}{\alpha})-C'_{\alpha}(\frac{\rho_n}{\alpha}))=0.
\end{equation*}
Using variable transformation gives $\int_{\alpha\rho_n}^{\alpha\rho_{n+1}}C_{\alpha}(x)\mathrm{d}x=0$. Similarly, the last identity can be proven. The proof is completed.

\textbf{Remark 3.4}. On the computation of the first positive zeros
of $C_{\alpha}(x)$ and $S_{\alpha}(x)$, we now don't have a good
method. Here, I propose an approximate approach to deal with it. Of course, it is not strict.
Denote $\eta_1(\alpha)$ as the first positive zero of
$C_{\alpha}(x)$ which can be considered as a function of $\alpha$.
Here, what we want to do is to find its first order approximation. Since
$C_{\alpha}(x)$ is the solution of the equation $y''(x)=-\alpha
y(\alpha^2x)$ with the initial conditions $y(0)=1$ and $ y'(0)=0$, we know
that if $\alpha=0$, the solution is $y(x)=1$ which has no zero, and
if $\alpha=1$, the solution is $y(x)=\cos(x)$ whose first positive
zero is $\frac{\pi}{2}$. Therefore, we can conclude that
$\eta_1(\alpha)$ is not an analytic function of $\alpha$, and has a
singularity at $\alpha=0$. Based on this observation, we assume that
$\eta_1(\alpha)$ has the following form
\begin{equation}
\eta_1(\alpha)=\alpha^{\kappa}\sum_{m=0}^{+\infty}x_m\alpha^m,
\end{equation}
where $\kappa$ is the index of singularity and $x_m's$ are the parameters undetermined. Under the first approximation, we can assume that $x_0+x_1=\frac{\pi}{2}$. Substituting $\eta_1(\alpha)$ into the
power series expansion of $C_{\alpha}(x)$ yields
\begin{equation*}
C_{\alpha}(\eta_1)=1-\frac{\alpha^{2\kappa+1}}{2}(x_0+x_1\alpha+\cdots)^2+\cdots=0.
\end{equation*}
By taking the first order approximation, we have
\begin{equation*}
1-\frac{\alpha^{2\kappa+1}}{2}x^2_0=0,
\end{equation*}
which means
\begin{equation*}
2\kappa+1=0,
\end{equation*}
and then
\begin{equation*}
1-\frac{x_0^2}{2}=0.
\end{equation*}
So we obtain $\kappa=-\frac{1}{2}$ and $x_0=\sqrt{2}$, and hence give the first order approximation $\eta_1(\alpha)=\alpha^{-\frac{1}{2}}(\sqrt2+x_1\alpha)$. By using $x_1=\frac{\pi}{2}-x_0$, we have the following interesting formula
\begin{equation}
\eta_1(\alpha)\simeq \alpha^{-\frac{1}{2}}\{(1-\alpha)\sqrt2+\alpha\frac{\pi}{2}\}.
\end{equation}
We should notice that $\sqrt2$ and $\frac{\pi}{2}$ are respectively just the first positive zeros of the solution  $y(x)=1-\frac{x^2}{2}$ (for $\alpha=0$) and the solution $y(x)=\cos(x)$ (for $\alpha=1$) of the equation
 $y''(x)=-y(\alpha x)$ with $y(0)=1,y'(0)=0$.

 By using the same method to deal with $\frac{S_{\alpha}(x)}{x}$,  we have the first order approximation formula of the first positive zero $\rho_1$ of $S_{\alpha}(x)$,
\begin{equation}
\rho_1(\alpha)\simeq \alpha^{-\frac{3}{2}}\{(1-\alpha)\sqrt6+\alpha\pi\}.
\end{equation}

Next we consider the number theoretical
properties of the zeroes of the sine-like function $S_{\alpha}(x)$.

\textbf{Proposition 3.6}.
 Notice that $\rho_0=0,\pm\rho_1,\cdots, \pm\rho_n, \cdots$ are all real zeroes of the function $S_{\alpha}(x)$. We will have
the following interesting formula
\begin{equation}
\frac{1}{\rho^2_1}+\frac{1}{\rho^2_2}+\cdots+\frac{1}{\rho^2_n}+\cdots=\frac{\alpha^3}{6},
\end{equation}
\begin{equation}
\frac{1}{\rho^{4}_1}+\frac{1}{\rho^{4}_2}+\cdots+\frac{1}{\rho^{4}_n}+\cdots=\frac{1}{36}\alpha^{6}
-\frac{1}{60}\alpha^{10}.
\end{equation}

\textbf{Proof}. Since the order of $S_{\alpha}(z)$ is zero, we know that for any positive real number $\epsilon$,
the following series
\begin{equation*}
\frac{1}{\rho^{\epsilon}_1}+\frac{1}{\rho^{\epsilon}_2}+\cdots+\frac{1}{\rho^{\epsilon}_n}+\cdots
\end{equation*}
is convergent.  Then by Hadamard's factorization theorem,  the infinite product of $S_{\alpha}(x)$ can be written as
\begin{equation*}
\frac{S_{\alpha}(x)}{x}=(1-\frac{x}{\rho_1})(1+\frac{x}{\rho_1})\cdots(1-\frac{x}{\rho_n})(1+\frac{x}{\rho_n})\cdots,
\end{equation*}
and by the series expansion of $S_{\alpha}(x)$, we can
obtain the formula (41). By considering the fourth order terms, we can get the formula (42). The proof is completed.

\textbf{Remark 3.5}. It is easy to see that this formula gives the classical
Euler's formulas if we take $\alpha=1$,
\begin{equation*}
1+\frac{1}{2^2}+\frac{1}{3^2}+\cdots=\frac{\pi^2}{6},
\end{equation*}
 and
\begin{equation*}
1+\frac{1}{2^4}+\frac{1}{3^4}+\cdots=\frac{\pi^4}{90}.
\end{equation*}

Of course, by some efforts, we can give high order formulas. If  more detailed information of zero
points is known, we can get more exact formulas.  I hope
that these secrets will be recovered in future.

\textbf{Open problem 3.1}. Give the exact values of the zeros of the
functions $S_{\alpha}(x)$ and $C_{\alpha}(x)$.
 
 In another paper [34], I will give the asymptotic formulas of zeros of $S_{\alpha}(x)$ and $C_{\alpha}(x)$.

\section{Power series method for the first order linear non-homogenous pantograph equation }

If non-homogenous term is a general smooth function, we do not know
how to solve the first order  linear pantograph equation since the
method of the variation of constant can not be applied in this case. Indeed, for the non-homogenous pantograph equation $y'(x)=\beta y(\alpha x)+q(x)$, the general solution of corresponding homogenous equation is $y(x)=cE_{\alpha}(\beta x)$. According to the variation of constant method, we can suppose that the special solution of the non-homogenous equation has the form $y(x)=c(x)E_{\alpha}(\beta x)$, and hence we have $c'(x)E_{\alpha}(\beta x)+c(x)\beta E_{\alpha}(\alpha\beta x)=c(\alpha x)\beta E_{\alpha}(\alpha\beta x)+q(x)$ form which we cannot eliminate two terms $c(x)\beta E_{\alpha}(\alpha\beta x)$ and $c(\alpha x)\beta E_{\alpha}(\alpha\beta x)$. Thus, the variation of constant method is invalid for pantograph equation, and then we must try to find other method.
In this section, we try to use power series to solve the first order linear
non-homogenous pantograph equation  with analytic non-homogenous
term. The power series method is simple but so powerful that we
use it to find out luckily the structure of the solutions for some important linear pantograph
equations.

\textbf{Theorem 4.1}. Consider the following pantograph equation
\begin{equation}
y'(x)=\beta y(\alpha x)+q(x), y(0)=a_0,
\end{equation}
where $\beta$ is a constant, and $q(x)$ is a analytic function with expansion
\begin{equation*}
q(x)=\sum_{n=0}^{+\infty}q_nx^n.
\end{equation*}
Then its solution is given by
 \begin{equation}
y(x)=a_0E_{\alpha}(\beta x)+\sum_{n=1}^{+\infty}\{\sum_{k=0}^{n-1}\frac{k!
\beta^{n-k-1}}{\alpha^{\frac{k(k+1)}{2}}}q_k\}\frac{\alpha^{\frac{n(n-1)}{2}}}{n!}x^n,
\end{equation}
or another form, if the double summations in (44) can exchange each other,
 \begin{equation}
y(x)=\{a_0+\sum_{k=0}^{+\infty}\frac{k!q_k}
{\beta^{k+1}\alpha^{\frac{k(k+1)}{2}}}\}E_{\alpha}(\beta x)- \sum_{k=0}^{+\infty}\frac{k!q_k}
{\beta^{k+1}\alpha^{\frac{k(k+1)}{2}}}\sum_{n=0}^{k}\frac{\alpha^{\frac{n(n-1)}{2}}}{n!}
\beta^nx^n.
\end{equation}

\textbf{Proof}. Assuming that the power series expansion of $y(x)$ is
\begin{equation*}
y(x)=\sum_{n=0}^{+\infty}a_nx^n,
\end{equation*}
and substituting it into the above equation and setting all coefficients of each term $x^n$ to be zeros,
we obtain the formulas of $a_n$ as follows
\begin{equation*}
a_{n+1}=\frac{\beta\alpha^n}{n+1}a_n+\frac{q_n}{n+1}.
\end{equation*}
Furthermore, we have (for $n\geq 1$)
\begin{equation*}
a_{n}=\frac{\alpha^{\frac{n(n-1)}{2}}}{n!}\{a_0\beta^{n}+\sum_{k=0}^{n-1}\frac{k!
\beta^{n-k-1}}{\alpha^{\frac{k(k+1)}{2}}}q_k\}.
\end{equation*}
So the solution can be represented by
 \begin{equation*}
y(x)=\sum_{n=0}^{+\infty}a_0\frac{\alpha^{\frac{n(n-1)}{2}}}{n!}\beta^{n}x^n+\sum_{n=1}^{+\infty}\{\sum_{k=0}^{n-1}\frac{k!
\beta^{n-k-1}}{\alpha^{\frac{k(k+1)}{2}}}q_k\}\frac{\alpha^{\frac{n(n-1)}{2}}}{n!}x^n
\end{equation*}
\begin{equation*}
=a_0E_{\alpha}(\beta x)+\sum_{n=1}^{+\infty}\{\sum_{k=0}^{n-1}\frac{k!
\beta^{n-k-1}}{\alpha^{\frac{k(k+1)}{2}}}q_k\}\frac{\alpha^{\frac{n(n-1)}{2}}}{n!}x^n,
\end{equation*}
or by exchanging the summations order
\begin{equation*}
y(x)=\{a_0+\sum_{k=0}^{+\infty}\frac{k!q_k}
{\beta^{k+1}\alpha^{\frac{k(k+1)}{2}}}\}E_{\alpha}(\beta x)- \sum_{k=0}^{+\infty}\frac{k!q_k}
{\beta^{k+1}\alpha^{\frac{k(k+1)}{2}}}\sum_{n=0}^{k}\frac{\alpha^{\frac{n(n-1)}{2}}}{n!}
\beta^nx^n.
\end{equation*}
The proof is completed.

 \textbf{Corollary 4.1}. If $q(x)$ is a polynomial of $n$ degree, the last summation will include only finite terms, that is,
\begin{equation}
y(x)=\{a_0+\sum_{k=0}^{n}\frac{k!q_k}
{\beta^{k+1}\alpha^{\frac{k(k+1)}{2}}}\}E_{\alpha}(\beta x)- \sum_{k=0}^{n}\frac{k!q_k}
{\beta^{k+1}\alpha^{\frac{k(k+1)}{2}}}\sum_{n=0}^{k}\frac{\alpha^{\frac{n(n-1)}{2}}}{n!}
\beta^nx^n.
\end{equation}

 \textbf{Theorem 4.2}. If $q(x)=q E_{\alpha}(\beta\alpha x)$, the solution of Eq.(43) is given by
 \begin{equation}
y(x)=cE_{\alpha}(\beta x)+qxE_{\alpha}(\alpha\beta x),
\end{equation}
where $c=a_0$.

\textbf{Proof}. We can prove the result by using the theorem 4.1. In fact, according to the power series expansion of $E_{\alpha}(\alpha\beta x)$, for $q(x)=q E_{\alpha}(\beta\alpha x)$, we have
\begin{equation*}
q(x)=\sum_{k=0}^{+\infty}q_kx^k=q\sum_{k=0}^{+\infty}\frac{\alpha^{\frac{k(k+1)}{2}}}{k!}
\beta^kx^k
\end{equation*}
which gives
\begin{equation*}
q_k=q\frac{\alpha^{\frac{k(k+1)}{2}}}{k!}\beta^k.
\end{equation*}
Substituting it into the solution (44) yields
\begin{equation*}
y(x)=a_0E_{\alpha}(\beta x)+q\sum_{n=1}^{+\infty}\{\sum_{k=0}^{n-1}\beta^{n-1}\}\frac{\alpha^{\frac{n(n-1)}{2}}}{n!}x^n,
\end{equation*}
\begin{equation*}
=a_0E_{\alpha}(\beta x)+q\sum_{n=1}^{+\infty}\beta^{n-1}\frac{\alpha^{\frac{n(n-1)}{2}}}{(n-1)!}x^n,
\end{equation*}
\begin{equation*}
=a_0E_{\alpha}(\beta x)+qx\sum_{n=0}^{+\infty}\beta^{n}\frac{\alpha^{\frac{n(n+1)}{2}}}{n!}x^n,
\end{equation*}
\begin{equation*}
=a_0E_{\alpha}(\beta x)+qx\sum_{n=0}^{+\infty}\beta^{n}\alpha^n\frac{\alpha^{\frac{n(n-1)}{2}}}{n!}x^n,
\end{equation*}
\begin{equation*}
=a_0E_{\alpha}(\beta x)+qxE_{\alpha}(\alpha\beta x).
\end{equation*}
The proof is completed.

\textbf{Remark 4.1.} Here we can give another derivation by
undetermined coefficients method. Assume the solution has the form
\begin{equation*}
y(x)=cE_{\alpha}(\gamma_1x)+BxE_{\alpha}(\gamma_2 x),
\end{equation*}
where $c, B, \gamma_1$ and $\gamma_2$ are unknown parameters. Substituting it into equation and setting
 the coefficients of $E_{\alpha}(\alpha\gamma_1x)$ and $E_{\alpha}(\alpha\gamma_2x)$ to be zeros yields that $c$
 is an arbitrary constant and $\gamma_1=\beta, \gamma_2=\alpha\beta$, and then $B=q$.

\textbf{Theorem 4.3}. If $q(x)=qx^k E_{\alpha}(\beta\alpha^{k+1} x)$ where $k\neq -1$, then the solution of Eq.(43) is given by
 \begin{equation*}
y(x)=cE_{\alpha}(\beta x)+\frac{q}{k+1}x^{k+1}E_{\alpha}(\alpha^{k+1}\beta x),
\end{equation*}
where $c=a_0$.

\textbf{Proof}. Assume that a special solution of the Eq.(43) is
\begin{equation*}
y^{*}(x)=Ax^hE_{\alpha}(\gamma x).
\end{equation*}
Then, we have
\begin{equation*}
\frac{d}{dx}y^{*}(x)=Ahx^{h-1}E_{\alpha}(\gamma x)+Ax^h\gamma E_{\alpha}(\alpha\gamma x),
\end{equation*}
and
\begin{equation*}
 y^{*}(\alpha x)=A\alpha^hx^hE_{\alpha}(\alpha\gamma x).
\end{equation*}
Therefore, we get
\begin{equation*}
Ahx^{h-1}E_{\alpha}(\gamma x)+Ax^h\gamma E_{\alpha}(\alpha\gamma x)
=\beta A\alpha^hx^hE_{\alpha}(\alpha\gamma x)+qx^k E_{\alpha}(\beta\alpha^{k+1} x),
\end{equation*}
and hence we have
\begin{equation*}
\gamma=\beta \alpha^h, qx^k=hAx^{h-1}, \gamma=\beta\alpha^{k+1},
\end{equation*}
from which it follows that
\begin{equation*}
h=k+1, \gamma=\beta \alpha^{k+1}, A=\frac{q}{k+1}.
\end{equation*}
Therefore, the special solution is
\begin{equation*}
y^{*}(x)=\frac{q}{k+1}x^{k+1}E_{\alpha}(\alpha^{k+1}\beta x),
\end{equation*}
from which we get the conclusion. The proof is completed.

 \textbf{Remark 4.2}. In some special cases, for the first order variable coefficient
equations, the exact solutions can be obtained. For example,
assuming that the solution of equation $y'(x)=p(x)y(\alpha x)$ is
$y(x)=E_{\beta}(\phi(x))$,
 then we have
\begin{equation*}
y'(x)=p(x)E_{\beta}(\phi(\alpha x))=\phi'(x)E_{\beta}(\beta\phi(x)),
\end{equation*}
which gives
\begin{equation*}
\beta=\alpha^{\gamma}, \phi(x)=Ax^{\gamma}, p(x)=\gamma
Ax^{\gamma-1},
\end{equation*}
that is, $y(x)=y(0)E_{\alpha^{\gamma}}(Ax^{\gamma})$ is the solution of
equation $y'(x)=\gamma Ax^{\gamma-1}y(\alpha x)$.

\section{The initial value problem of the second order  linear pantograph equation}
For the second order linear pantograph equation without
the first order term
\begin{equation}
y''(x)= A y(\gamma x),
\end{equation}
\begin{equation}
y(0)=c_1, y'(0)=c_2,
\end{equation}
where $0<\gamma<1$, we can solve it by the special function $E_{\alpha}(x)$. In fact, taking $y(x)=E_{\alpha}(\beta x)$ and using
 \begin{equation*}
y''(x)= \alpha\beta^2 y(\alpha^2 x),
\end{equation*}
gives
\begin{equation*}
\alpha\beta^2=A,
\end{equation*}
\begin{equation*}
\alpha^2=\gamma.
\end{equation*}
Solving them yields $\alpha=\sqrt{\gamma}$ and $\beta=\pm\sqrt{\frac{A}{\sqrt{\gamma}}}$. Therefore, we obtain
two basic solutions
\begin{equation*}
y_1(x)=E_{\sqrt{\gamma}}(\sqrt{\frac{A}{\sqrt{\gamma}}} x),
\end{equation*}
\begin{equation*}
y_2(x)=E_{\sqrt{\gamma}}(-\sqrt{\frac{A}{\sqrt{\gamma}}} x),
\end{equation*}
and hence the general solution can be represented by
\begin{equation}
y(x)=a_1E_{\sqrt{\gamma}}(\sqrt{\frac{A}{\sqrt{\gamma}}} x)+a_2E_{\sqrt{\gamma}}(-\sqrt{\frac{A}{\sqrt{\gamma}}} x),
\end{equation}
where $a_1=\frac{1}{2}(c_1+c_2\sqrt{\frac{\sqrt\gamma}{A}})$ and $a_2=\frac{1}{2}(c_1-c_2\sqrt{\frac{\sqrt\gamma}{A}})$.

We now consider the general case and give the following results.

\textbf{Theorem 5.1.} For the second order linear pantograph  equation with the first order term,
\begin{equation}
y''(x)+py'(\alpha x)+ qy(\alpha^2 x)=0,
\end{equation}
\begin{equation}
y(0)=c_1, y'(0)=c_2,
\end{equation}
where $y'(\alpha x)=y'(t)|_{t=\alpha x}$. There are the following
two cases:

(i). If $\vartriangle=p^2-4\alpha q\neq0$, the solution is given by
\begin{equation}
y(x)=a_1E_{\alpha}(\beta_1 x)+a_2E_{\alpha}(\beta_2 x),
\end{equation}
where $\beta_{1,2}=\frac{-p\pm\sqrt{\vartriangle}}{2\alpha}, a_1=\frac{c_1\beta_2-c_2}{\beta_2-\beta_1}$
 and $a_2=\frac{c_1\beta_1-c_2}{\beta_1-\beta_2}$.

(ii). If $\vartriangle=p^2-4\alpha q=0$, the solution is given by
\begin{equation}
y(x)=a_1E_{\alpha}(\beta x)+a_2xE_{\alpha}(\alpha\beta x),
\end{equation}
where $\beta=-\frac{p}{2\alpha}, a_1=c_1$ and
$a_2=c_2-c_1\beta$.

\textbf{Proof}. Taking $y(x)=E_{\alpha}(\beta x)$ and substituting it into the above equation and eliminating
$E_{\alpha}(\beta\alpha^2 x)$
gives the characteristic equation
\begin{equation}
\alpha\beta^2+p\beta+q=0.
\end{equation}
Solving it yields
\begin{equation*}
\beta_{1,2}=\frac{-p\pm\sqrt{\vartriangle}}{2\alpha},
\end{equation*}
where $\vartriangle=p^2-4\alpha q$. Therefore, if $\vartriangle\neq 0$, we obtain two basic solutions
\begin{equation*}
y_1(x)=E_{\alpha}(\beta_1 x),
\end{equation*}
\begin{equation*}
y_2(x)=E_{\alpha}(\beta_2 x),
\end{equation*}
and then the general solution can be represented by
\begin{equation*}
y(x)=a_1E_{\alpha}(\beta_1 x)+a_2E_{\alpha}(\beta_2 x),
\end{equation*}
where the coefficients $a_1$ and $a_2$ can be determined by initial conditions.

If $\vartriangle=0$, we can give a basic solution, that is,
\begin{equation*}
y_1(x)=E_{\alpha}(\beta x).
\end{equation*}
In order to obtain another basic solution, since the usual method of
variation of constant is not suitable, we must find other method. In
section 7, we will use a kind of operator technics to deal with this
problem and
 give the details of the corresponding theory. Here, we only verify that another basic solution is
\begin{equation*}
y_2(x)=xE_{\alpha}(\alpha\beta x).
\end{equation*}
In fact, by inserting it into equation, it is easy to see that it is just a solution. Therefore, in this case,
 the general solution (54) is given. The proof is completed.

Next we consider the second order linear pantograph equations with some non-homogenous terms and give the following result.

\textbf{Theorem 5.2}. For the non-homogenous equation
\begin{equation*}
y''(x)+py'(\alpha x)+ qy(\alpha^2 x)=\sum_{k=1}^mA_kE_{\alpha}(r_kx),
\end{equation*}
where $A_k$ and $r_k$ are known constants, a special solution is given by
\begin{equation}
y^{*}(x)=\sum_{k=1}^m\frac{A_k\alpha^3}{r_k^2+pr_k\alpha+q\alpha^3}E_{\alpha}(\frac{r_k}{\alpha^2}x),
\end{equation}
where we require $r_k^2+pr_k\alpha+q\alpha^3\neq0$. Further, the general solution is the summation of the general solution of the homogenous equation and the special solution.

\textbf{Proof}. Assuming that $y_k^{*}(x)=B_kE_{\alpha}(s_kx)$ is the special solutions of the equation
$y''(x)+py'(\alpha x)+ qy(\alpha^2 x)=A_kE_{\alpha}(r_kx)$ and substituting it into the equation yields
\begin{equation*}
s_k=\frac{r_k}{\alpha^2}, B_k=\frac{A_k\alpha^3}{r_k^2+pr_k\alpha+q\alpha^3}.
\end{equation*}
Therefore, from the theorem 2.4, we get the special solution (56). The proof is completed.

It is well-known that there is a simple conservation law for the
usual vibration equation ($\alpha=1$)namely energy conservation law
according to the invariance of time translation. But we do not know whether
there exists a similar simple conservation law for the case of
$0<\alpha<1$. A complicated conservation law can be given by using
the addition formulas of $S_{\alpha}(x)$ and $C_{\alpha}(x)$.

\textbf{Theorem 5.3.}  For the second order pantograph equation
\begin{equation}
q''(t)= -\alpha q(\alpha^2 t),
\end{equation}
\begin{equation}
q(0)=q_0, q'(0)=v_0,
\end{equation}
if denote $p(t)=q'(\frac{t}{\alpha})$, then there is a conservation law as follows
\begin{equation}
I(t)=\frac{1}{q^2_0+v^2_0}\sum_{n=0}^{+\infty} \frac{(-1)^nt^{n}}{n!}(q_0 q^{(n)}(t)+v_0p^{(n)}(t))=1.
\end{equation}

 \textbf{Proof}. In the addition formula
 \begin{equation*}
C_{\alpha}(x+y)= \sum_{n=0}^{+\infty} \alpha^{n(2n-1)}\frac{(-1)^nx^{2n}}{(2n)!}C_{\alpha}
(\alpha^{2n}y)-\sum_{n=0}^{+\infty} \alpha^{n(2n+1)}\frac{(-1)^nx^{2n+1}}{(2n+1)!}S_{\alpha}(\alpha^{2n+1}y),
\end{equation*}
we take $y=-x=t$, and notice that $C^{(2n)}_{\alpha}(y)=(-1)^n\alpha^{n(2n-1)}C_{\alpha}(\alpha^{2n}y)$
and $C^{(2n+1)}_{\alpha}(y)=(-1)^{n+1}\alpha^{n(2n+1)}S_{\alpha}(\alpha^{2n+1}y)$,  then we have
 \begin{equation}
1=C_{\alpha}(t-t)= \sum_{n=0}^{+\infty} \frac{(-1)^nt^{n}}{n!}C^{(n)}_{\alpha}(t).
\end{equation}
Since the solution of the initial value problem is
\begin{equation*}
q(t)=q_0C_{\alpha}(t)+v_0S_{\alpha}(t),
\end{equation*}
we have
\begin{equation*}
C_{\alpha}(t)=\frac{q_0 q(t)+v_0q'(\frac{t}{\alpha})}{q^2_0+v^2_0}=\frac{q_0 q(t)+v_0p'(t)}{q^2_0+v^2_0}.
\end{equation*}
Then by  substituting it into Eq.(60) gives the result. The proof is completed.

\textbf{Remark 5.1}. For the general second order equation
$y''(x)+py'(\alpha x)+qy(\beta x)=0$ with $\beta\neq \alpha^2$, we
need introduce new special functions to solve it.

\textbf{Remark 5.2}. We almost don't consider the real solutions when some eigenvalues are perhaps not real numbers. By means of the relations between $E_{\alpha}(x)$ and $S_{\alpha}(x)$ and $C_{\alpha}(x)$, we can easily obtain the corresponding results. In the paper, we omit them for simplicity.

\textbf{Open problem 5.1}. Give a variational principle for the
second order  vibration equation $y''(x)=-k y(\alpha x)$.

\section{The system of the linear pantograph equations}

\textbf{Theorem 6.1}. For the system of linear pantograph equations
\begin{equation}
x'_n(t)=\beta x_n(\alpha t)+x_{n+1}(\alpha t), n=1,\cdots,m,
\end{equation}
where $x_{m+1}(t)=0$, its solutions of initial value problem at the original point $t=0$ are given by
\begin{equation}
x_{m-k}(t)=\sum_{j=0}^{k}\frac{x_{m-k+j}(0)}{j!}\alpha^{\frac{j(j-1)}{2}}t^jE_{\alpha}(\alpha^j\beta t),
\end{equation}
where $k=0,1,\cdots,m-1$.

\textbf{Proof}. Firstly, we prove the theorem for $k=0$. The
corresponding equations become
\begin{equation*}
x_m'(t)=\beta x_m(\alpha t),
\end{equation*}
whose solution is given by
\begin{equation*}
x_m(t)= x_m(0)E_{\alpha}(\beta t).
\end{equation*}
Further, we consider
\begin{equation*}
x_{m-1}'(t)=\beta x_{m-1}(\alpha t)+x_m(\alpha t),
\end{equation*}
that is
\begin{equation*}
x_{m-1}'(t)=\beta x_{m-1}(\alpha t)+x_m(0)E_{\alpha}(\beta t).
\end{equation*}
By the theorem 4.2, we have
\begin{equation*}
x_{m-1}(t)= x_{m-1}(0)E_{\alpha}(\beta t)+x_m(0)tE_{\alpha}(\alpha\beta t).
\end{equation*}
In general, we assume that the formula (62) holds for $k$, then we prove it holds for $k+1$. In fact, we have
\begin{equation*}
x_{m-k-1}'(t)=\beta x_{m-k-1}(\alpha t)+x_{m-k}(\alpha t),
\end{equation*}
that is,
\begin{equation*}
x_{m-k-1}'(t)=\beta x_{m-k-1}(\alpha t)+\sum_{j=0}^{k}\frac{x_{m-k+j}(0)}{j!}\alpha^{\frac{j(j+1)}{2}}t^jE_{\alpha}(\alpha^{j+1}\beta t).
\end{equation*}
 From the theorems 2.4 and 4.3, we have
\begin{equation*}
x_{m-k-1}(t)=x_{m-k-1}(0)E_{\alpha}(\beta t)+\sum_{j=0}^{k}\frac{x_{m-k+j}(0)}{j!}\alpha^{\frac{j(j+1)}{2}}\frac{1}{j+1}t^{j+1}E_{\alpha}(\alpha^{j+1}\beta t)
\end{equation*}
\begin{equation*}
=x_{m-k-1}(0)E_{\alpha}(\beta t)+\sum_{j=0}^{k}\frac{x_{m-k+j}(0)}{(j+1)!}\alpha^{\frac{j(j+1)}{2}}t^{j+1}E_{\alpha}(\alpha^{j+1}\beta t)
\end{equation*}
\begin{equation*}
=\sum_{j=0}^{k+1}\frac{x_{m-k+j-1}(0)}{j!}\alpha^{\frac{j(j-1)}{2}}t^jE_{\alpha}(\alpha^j\beta t).
\end{equation*}
By the mathematical induction method, the proof is completed.

According to the theorem 6.1, we have the following theorem.

\textbf{Theorem 6.2}. Consider the system of linear pantograph equations
\begin{equation}
\frac{\mathrm{d}Y(t)}{\mathrm{d}t}=AY(\alpha t),
\end{equation}
with initial condition
\begin{equation}
Y(0)=Y_0,
\end{equation}
where $y(t)=(y_1(t),\cdots,y_n(t))^T$ is a vector function and $A=(a_{ij})_{n\times n}$ is a constant matrix.

\textbf{Case (i)}. If $A$ has $n$ distinct  eigenvalues $\lambda_1,\cdots,\lambda_n$, then there is an invertible matrix $P$ such that $P^{-1}AP$ is diagonal, i.e.,
\begin{equation*}
P^{-1}AP=\left(
  \begin{array}{ccc}
    \lambda_1  &\cdots &0 \\
   \vdots &\ddots&\vdots\\
   0&\cdots & \lambda_n\\
  \end{array}
\right).
\end{equation*}
and hence the equations become
\begin{equation*}
x_j'(t)=\lambda_jx_j(\alpha t), j=1,\cdots,n,
\end{equation*}
whose solutions are given by
\begin{equation}
x_j(t)=x_j(0)E_{\alpha}(\lambda_j t), j=1,\cdots,n.
\end{equation}
Respectively, we have the solutions $Y(t)=PX(t)$, and the general solution of the homogenous equation
is the linear combination of these solutions.

\textbf{Case (ii)}. If $A$ has $m$ distinct eigenvalues are $\lambda_1,\cdots,\lambda_r$ with multiplies $n_1,\cdots,n_r$, and corresponding elementary factors are
 $(\lambda-\lambda_1)^{k_{11}}, \cdots, (\lambda-\lambda_1)^{k_{1m_1}},\cdots,(\lambda-\lambda_r)^{k_{r1}}, \cdots, (\lambda-\lambda_r)^{k_{rm_r}}$ where $k_{j1}+...+k_{jm_j}=n_j$ for $j=1,\cdots,r$, then there is an invertible
 matrix $P$ such that $P^{-1}AP$ is a Jordan form matrix, that is,
 \begin{equation*}
P^{-1}AP=\left(
  \begin{array}{ccc}
    J_1  &\cdots &0 \\
   \vdots &\ddots&\vdots\\
   0&\cdots & J_{m_1+\cdots+m_r}\\
  \end{array}
\right).
\end{equation*}
 where $J_k's$ are the corresponding Jordan blocks. For every Jordan block,  the solutions $X$ can be given by formula (62).
 Then  we have solutions $Y=PX$, and the general solution of the homogenous equation is the linear
 combination of these solutions.

\section{Operator method}

We use an operator method to deal with linear pantograph equations.
 Two basic operators are derivative operator $D$ and the scale operator $ T_{\alpha}$ defined by
\begin{equation*}
Dy(x)=\frac{\mathrm{d}y}{\mathrm{d}x},
\end{equation*}
\begin{equation*}
T_{\alpha}y(x)=y(\alpha x).
\end{equation*}
Two basic properties of these two operators are
\begin{equation*}
DT_{\alpha}=\alpha T_{\alpha} D,
\end{equation*}
\begin{equation*}
T_{\alpha_1}T_{\alpha_2}=T_{\alpha_1\alpha_2}.
\end{equation*}

\textbf{Theorem 7.1}. For a special  first order equation
\begin{equation}
Dy(x)=pT_{\alpha}y(x)+qE_{\alpha}(\gamma x),
\end{equation}
the general solution is given by
\begin{equation}
y(x)=cE_{\alpha}(px)+\frac{q\alpha}{\gamma-p\alpha}E_{\alpha}(\frac{\gamma}{\alpha}x),
\end{equation}
where $\gamma\neq p\alpha$ and $c$ is an arbitrary constant which can be determined by $y(0)$.

\textbf{Proof}. Formally, we have
\begin{equation*}
y(x)=q(D-pT_{\alpha})^{-1}E_{\alpha}(\gamma x).
\end{equation*}
To compute the right side of the above equation, we use the following formula
\begin{equation*}
(D-pT_{\alpha})E_{\alpha}(\beta x)=(\beta-p)E_{\alpha}(\alpha\beta x).
\end{equation*}
So, taking $\beta=\frac{\gamma}{\alpha}$ gives
\begin{equation*}
(D-pT_{\alpha})^{-1}E_{\alpha}(\gamma x)=\frac{1}{\frac{\gamma}{\alpha}-p}E_{\alpha}(\frac{\gamma}{\alpha}x).
\end{equation*}
Therefore, a special solution is given by
\begin{equation*}
y^*(x)=\frac{q\alpha}{\gamma-p\alpha}E_{\alpha}(\frac{\gamma}{\alpha}x),
\end{equation*}
and the general solution is
\begin{equation*}
y(x)=cE_{\alpha}(px)+\frac{q\alpha}{\gamma-p\alpha}E_{\alpha}(\frac{\gamma}{\alpha}x),
\end{equation*}
where $c$ is an arbitrary constant. The proof is completed.

\textbf{Remark 7.1}. For $\gamma= p\alpha$, by theorem 4.2, the solution is given by
\begin{equation*}
y(x)=cE_{\alpha}(px)+qxE_{\alpha}(p\alpha x).
\end{equation*}

Now we consider the general second order linear homogenous pantograph equation
 \begin{equation}
y''(x)+py'(\alpha x)+ qy(\alpha^2 x)=0,
\end{equation}
whose operator form is
 \begin{equation*}
(D^2+pT_{\alpha}D+qT^2_{\alpha})y(x)=0.
\end{equation*}
In general, we have the operator decomposition
\begin{equation*}
D^2+pT_{\alpha}D+qT^2_{\alpha}=(D-\lambda_1T_{\alpha})(D-\lambda_2T_{\alpha}),
\end{equation*}
where the parameters $\lambda_1$ and  $\lambda_1$ satisfy
\begin{equation*}
\lambda_1+\lambda_2\alpha=-p, \lambda_1\lambda_2=q.
\end{equation*}

\textbf{Remark 7.2}. If we take $\lambda_1=\alpha\beta_1$ and $\lambda_2=\beta_2$, then we will give the same results with theorem 5.1.

The key of solving this equation is the following lemma which can be easily proven.

\textbf{Lemma 7.1}. If $y(x)$ is the solution of equation $(D-\lambda_2T_{\alpha})y(x)=0$, it is
also the solution of the equation $(D-\lambda_1T_{\alpha})(D-\lambda_2T_{\alpha})y(x)=0$. If
 $z(x)$ is the solution of equation $(D-\lambda_1T_{\alpha})z(x)=0$, and $y(x)$ is the solution
 of equation $(D-\lambda_2T_{\alpha})y(x)=z(x)$, then $y(x)$ is also the solution of the
 equation $(D-\lambda_1T_{\alpha})(D-\lambda_2T_{\alpha})y(x)=0$.

By the lemma 7.1, we have the following theorem.

\textbf{Theorem 7.2}. For the initial value problem at the original point $x=0$ of equation (68), there are the following two cases.

\emph{\textbf{Case (i)}}. If $\lambda_1\neq\lambda_2\alpha$, the solution is give by
\begin{equation*}
y(x)=\frac{c_1}{\frac{\lambda_1}{\alpha}-\lambda_2}E_{\alpha}(\frac{\lambda_1}{\alpha}x)+c_2E_{\alpha}(\lambda_2x),
\end{equation*}
where $c_1$ and $c_2$ can be determined by initial conditions.

\emph{\textbf{Case (ii)}}.  If $\lambda_1=\lambda_2\alpha$, the solution is give by
\begin{equation*}
y(x)=c_1 E_{\alpha}(\lambda_1 x)+c_2xE_{\alpha}(\alpha\lambda_1 x),
\end{equation*}
where $c_1$ and $c_2$ can be determined by initial conditions.

\textbf{Proof}.  In fact, by lemma 7.1, the first basic solution $y_1(x)$ satisfies
\begin{equation*}
(D-\lambda_2T_{\alpha})y_1(x)=0,
\end{equation*}
and then
\begin{equation*}
y_1(x)=c_1E_{\alpha}(\lambda_2x).
\end{equation*}
Therefore, the second basic solution satisfies
\begin{equation*}
(D-\lambda_2T_{\alpha})y(x)=c_2E_{\alpha}(\lambda_1x),
\end{equation*}
whose special solution is given by
\begin{equation*}
y^*(x)=c_2(D-\lambda_2T_{\alpha})^{-1}E_{\alpha}(\lambda_1x).
\end{equation*}
From formula (67), if $\lambda_1\neq\lambda_2\alpha$, we know
\begin{equation*}
y^*(x)=\frac{c_2}{\frac{\lambda_1}{\alpha}-\lambda_2}E_{\alpha}(\frac{\lambda_1}{\alpha}x).
\end{equation*}
 Hence, in this case, the general solution is given by
\begin{equation*}
y(x)=c_1E_{\alpha}(\lambda_2x)+\frac{c_2}{\frac{\lambda_1}{\alpha}-\lambda_2}E_{\alpha}(\frac{\lambda_1}{\alpha}x),
\end{equation*}
where $c_1$ and $c_2$ are two arbitrary constants which can be determined by $y(0)$ and $y'(0)$.

If $\lambda_1=\lambda_2\alpha$, from formula (67), $y(x)$ is meaningless, so this method is invalid.
In this case, we need to solve
\begin{equation}
(D-\lambda_1T_{\alpha})y(x)=E_{\alpha}(\lambda_1\alpha x).
\end{equation}
By the power series method in section 4, we can assume that the solution has the following form
 \begin{equation*}
y(x)=c_1E_{\alpha}(s_1x)+c_2xE_{\alpha}(s_2x),
\end{equation*}
and insert it into the above equation to obtain
\begin{equation*}
s_1=\lambda_2, s_2=\alpha\lambda_2=\lambda_1.
\end{equation*}
So we have
\begin{equation*}
y(x)=c_1 E_{\alpha}(\lambda_2 x)+c_2xE_{\alpha}(\alpha\lambda_2 x),
\end{equation*}
where $c_1$ and $c_2$ can be determined by $y(0)$ and $y'(0)$. The proof is completed.

In general, for the  $n$-th order  linear pantograph
equation
\begin{equation}
y^{(n)}(t)+p_{n-1}y^{(n-1)}(\alpha t)+p_{n-2}y^{(n-2)}(\alpha^2
t)+\cdots+p_1y'(\alpha^{n-1} t)+p_0y(\alpha^n t)=0,
\end{equation}
we can write it as follows
\begin{equation*}
D^ny(x)+p_{n-1}T_\alpha D^{n-1}y(x)+p_{n-2}T^2_\alpha D^{n-2}y(x)\cdots+p_1T_{\alpha}^{n-1}Dy(x)+p_0T^n_\alpha y(x)=0,
\end{equation*}
where $p_i's$ are constants for $i=0,\cdots, n-1$, and $T^k_{\alpha}=T_{\alpha}\cdots T_{\alpha}=T_{\alpha^k}$.
Furthermore, we have the following two ways to deal with it.

The first way is to  solve it directly.  By taking $y(x)=E_{\alpha}(\beta t)$ and substituting it into the
above equation and eliminating $E_{\alpha}(\beta \alpha^n t)$, we get the characteristic equation
\begin{equation}
\alpha^{\frac{n(n-1)}{2}}\beta^n+p_{n-1}\alpha^{\frac{(n-1)(n-2)}{2}}\beta^{n-1}+\cdots+p_1\beta +p_0=0.
\end{equation}
If this equation has no multiple roots, that is, it has $n$ distinct roots $\beta_j$ for $j=1,\cdots,n$, we
can give the general solution by
\begin{equation}
y(t)=a_1E_{\alpha}(\beta_1 t)+\cdots+a_nE_{\alpha}(\beta_n t),
\end{equation}
where $a_1,\cdots,a_n$ can be uniquely determined by the initial conditions at the original point. If there are some multiple
roots, we will use the corresponding operator method to deal with it. By the similar method, we can prove
the following general theorem.

\textbf{Theorem 7.3}. If $\beta_1,\cdots,\beta_m$ are the distinct roots of the characteristic equation
with multiplicities respectively $n_1,\cdots,n_m$ satisfying $n_1+\cdots+n_m=n$. Then the following $n$ functions
 make up the basic solutions system of the initial value problem at original point of the equation (70),
\begin{equation*}
y_1(t)=E_{\alpha}(\beta_1 t), y_2(t)=t E_{\alpha}(\beta_1\alpha t),\cdots, y_{n_1}(t)=t^{n_1-1}E_{\alpha}(\beta_1 \alpha^{n_1-1} t),
\end{equation*}
\begin{equation*}
y_{n_1+1}(t)=E_{\alpha}(\beta_2 t), \cdots, y_{n_1+n_2}(t)=t^{n_2-1}E_{\alpha}(\beta_2 \alpha^{n_2-1} t),
\end{equation*}
\begin{equation*}
\cdots
\end{equation*}
\begin{equation}
y_{n-n_m+1}(t)=E_{\alpha}(\beta_m t),\cdots, y_{n}(t)=t^{n_m-1}E_{\alpha}(\beta_m \alpha^{n_m-1} t).
\end{equation}
In other words, the solution of the initial value problem at original point of the equation (70) is the linear combination of these basic solutions with coefficients which can be determined by initial conditions.

\textbf{Proof}. We only take $n=3$ and suppose that $\lambda_0$ is a characteristic root with multiplicity three to prove the result. The general case can be proved by the similar method. Denote $Q_3(D, T_{\alpha})$ and $Q_3(\lambda)$ respectively as
\begin{equation*}
Q_3(D, T_{\alpha})=D^3+p_2T_{\alpha}D^2+p_1T_{\alpha^2}D+p_0T_{\alpha^3},
\end{equation*}
\begin{equation*}
Q_3(\lambda)=\alpha^3\lambda^3+p_2\alpha\lambda^2+p_1\lambda+p_0.
\end{equation*}
By direct computation, we have
\begin{equation*}
Q_3(D, T_{\alpha})(x^2E_{\alpha}(\lambda_0\alpha^2 x))=Q_3(\lambda_0)\alpha^6x^2E_{\alpha}(\lambda_0\alpha^5 x)
\end{equation*}
\begin{equation*}
+Q'_3(\lambda_0)2\alpha xE_{\alpha}(\lambda_0\alpha^4 x)+Q''_3(\lambda_0)E_{\alpha}(\lambda_0\alpha^3 x).
\end{equation*}
Furthermore, since $\lambda_0$ is a characteristic root with multiplicity three, we have
\begin{equation*}
Q_3(\lambda_0)=Q'_3(\lambda_0)=Q''_3(\lambda_0)=0,
\end{equation*}
and hence, it follows that
\begin{equation*}
Q_3(D, T_{\alpha})(x^2E_{\alpha}(\lambda_0\alpha^2 x))=0.
\end{equation*}
This means that $x^2E_{\alpha}(\lambda_0\alpha^2 x)$ is a basic solution. By the same reason, $xE_{\alpha}(\lambda_0\alpha x)$ and $E_{\alpha}(\lambda_0 x)$ are other two basic solutions. Therefore, the general solution can be given by
\begin{equation}
y(x)=c_1E_{\alpha}(\lambda_0 x)+c_2xE_{\alpha}(\lambda_0\alpha x)+c_3x^2E_{\alpha}(\lambda_0\alpha^2 x),
\end{equation}
where $c_1,c_2$ and $c_3$ can be determined by the initial conditions. The proof is completed.

The second way is to transform Eq.(70) into the system of linear pantograph equations. In fact,
letting $x_1(t)=y(t), x_2(t)=x'_1(\frac{t}{\alpha}),\cdots, x_n(t)=x'_{n-1}(\frac{t}{\alpha})$, then the
equation (70) is transformed into the following  system of the first order linear equations
\begin{equation*}
x'_1(t)=x_2(\alpha t),
\end{equation*}
\begin{equation*}
x'_2(t)=x_3(\alpha t),
\end{equation*}
\begin{equation*}
\cdots
\end{equation*}
\begin{equation}
x'_n(t)=-\frac{p_{0}}{\alpha^{n-1}}x_{1}(\alpha t)-\frac{p_{1}}{\alpha^{n-1}}x_{2}(\alpha t)
-\cdots-\frac{p_{n-1}}{\alpha^{n-1}}x_{n}(\alpha t).
\end{equation}
For example, the second order pantograph  equation
\begin{equation*}
q''(t)=-kq(\alpha^2 t),
\end{equation*}
can be transformed to the Hamiltonian-like form
\begin{equation*}
q'(t)=p(\alpha t),
\end{equation*}
\begin{equation*}
p'(t)=-\frac{k}{\alpha}q(\alpha t).
\end{equation*}
This means that for the initial value problem at original point $t=0$, the theory of the high order linear pantograph equations is equivalent to the theory
of the system of the first order linear pantograph equations. Therefore, we can give the solutions of the high
order pantograph equation by the solution of the system of the first order pantograph equations.

For simplicity, we denote
\begin{equation*}
P(D,T_\alpha)=D^n+p_{n-1}T_\alpha D^{n-1}+p_{n-2}T^2_\alpha D^{n-2}\cdots+p_1T_{\alpha}^{n-1}D+p_0T^n_\alpha.
\end{equation*}
We have the following result on the non-homogenous equation.

\textbf{Theorem 7.4}. For non-homogenous equation
\begin{equation}
P(D,T_\alpha)y(x)=\sum_{k=1}^mA_kE_{\alpha}(r_kx),
\end{equation}
where $A_k$ and $r_k$ are known constants, a special solution is given by
\begin{equation}
y^{*}(x)=\sum_{k=1}^mB_kE_{\alpha}(\frac{r_k}{\alpha^n}x),
\end{equation}
where
\begin{equation}
B_k=\frac{A_k}{\sum_{j=0}^np_{j}r^{j}_k\alpha^{-\frac{j(2n-j+1)}{2}}},
\end{equation}
where $p_n=1$ and we require
\begin{equation}
\sum_{j=0}^np_{j}r^{j}_k\alpha^{-\frac{j(2n-j+1)}{2}}\neq 0.
\end{equation}
Further, the general solution is given by the summation of the general solution of the homogenous equation and the special solution.

\textbf{Proof}. Assume that the special solution of the equation $P(D,T_\alpha)y(x)=A_kE_{\alpha}(r_kx)$ is $y_k(x)=B_kE_{\alpha}(s_kx)$. Then, we have
\begin{equation*}
y_k^{(j)}(x)=B_ks_k^l\alpha^{\frac{j(j-1)}{2}}E_{\alpha}(\alpha^js_kx),
\end{equation*}
and substitute it into the equation to get $s_k=\frac{r_k}{\alpha^n}$ and the values of $B_k$. The proof is completed.

When we take $n=2$, we get the theorem 5.2.

\textbf{Definition 7.1}. If $r_k$ satisfies the equation (79), it is called the resonance frequency, correspondingly, $A_kE_{\alpha}(r_kx)$ is called the resonance term.

Now we consider the non-homogenous linear pantograph equations with resonance terms. We only give the result on the second order equation. Other high order equation can be dealt with by the same method.

\textbf{Theorem 7.5}. Denote $Q(r)=r^2+pr\alpha+q\alpha^3$. For the non-homogenous pantograph equation
\begin{equation*}
y''(x)+py'(\alpha x)+qy(\alpha^2x)=AE_{\alpha}(rx),
\end{equation*}
where $Q(r)\neq0$ which means that $AE_{\alpha}(rx)$ is a resonance term. There are the following two cases:

\textbf{Case (i)}. $Q(r)=0$ and $Q'(r)=2r+p\alpha\neq0$. Then a special solution is given by
\begin{equation}
y(x)=\frac{A\alpha}{2r+p\alpha}xE_{\alpha}(\frac{r}{\alpha}x).
\end{equation}

\textbf{Case (ii)}. $Q(r)=0$ and $Q'(r)=2r+p\alpha=0$. Then a special solution is given by
\begin{equation}
y(x)=\frac{A}{2}x^2E_{\alpha}(rx).
\end{equation}

\textbf{Proof}. Firstly, we prove the result in case (i). Assume that the special solution has the form
\begin{equation*}
y(x)=BxE_{\alpha}(sx).
\end{equation*}
Then we have
\begin{equation*}
y'(x)=BE_{\alpha}(sx)+BsxE_{\alpha}(s\alpha x),
\end{equation*}
\begin{equation*}
y''(x)=2BsE_{\alpha}(s\alpha x)+Bs^2\alpha xE_{\alpha}(s\alpha^2 x).
\end{equation*}
Substituting these terms into the equation give $Q(r)=0$ and
\begin{equation*}
s=\frac{r}{\alpha}, B=\frac{A\alpha}{2r+p\alpha}.
\end{equation*}

Next we prove the case (ii). Assume that the special solution has the form
\begin{equation*}
y(x)=Bx^2E_{\alpha}(sx).
\end{equation*}
Then we have
\begin{equation*}
y'(x)=2BxE_{\alpha}(sx)+Bsx^2E_{\alpha}(s\alpha x),
\end{equation*}
\begin{equation*}
y''(x)=2BE_{\alpha}(sx)+4BsxE_{\alpha}(s\alpha x)+Bs^2\alpha x^2E_{\alpha}(s\alpha^2 x).
\end{equation*}
Substituting these terms into the equation give $Q(r)=0$ and $Q'(r)=0$ and
\begin{equation*}
s=r, B=\frac{A}{2}.
\end{equation*}
The proof is completed.

\textbf{Remark 7.1}. The theorems 7.4 and 7.5 only mean that if the solution exists, it will be given by the summation of the special solution and the general solution of the homogenous equation, but it doesn't mean that the solution of the non-homogenous equation (76) must exist. The existence, uniqueness and non-uniqueness of the non-homogenous pantograph equations will be discussed in section 8.

\section{Existence, uniqueness, non-uniqueness and representation of solutions at a general initial point}

In previous sections, all results are obtained under the initial conditions at the original point. If the initial condition is taken at a general point which is not the original point, whether do the existence and uniqueness of solution hold?
How to represent these solutions if they exist? How many solutions will exist if uniqueness does not hold? These problems are not trivial. In the section, we will give the results about these problems. For the purpose, we need an important theorem about the zeroes of $E_{\alpha}(x)$ which is perhaps considered as the most important property of the function which is equivalent to the theorem 6 in [8].

\textbf{Morris-Feldstein-Bowen-Hahn Theorem: first form[8].}  Every nontrivial solution of the equation $y'(x)=-y(\alpha x)$ has an infinity of positive zeroes.

Equivalently, in other words,  we write it as the following form by using the special function  $E_{\alpha}(x)$.

\textbf{Morris-Feldstein-Bowen-Hahn Theorem: second form.}(MFBH theorem for simplicity). $E_{\alpha}(x)$ has an infinity of negative zeros.

\textbf{Remark 8.1}. The above theorem is obtained and proven in [8] by Morris-Feldstein-Bowen, and an elementary proof is given by Hahn [8].

By the MFBH theorem, we can get the following theorem.

\textbf{Theorem 8.1}. For the initial value problem of linear equation
\begin{equation}
y'(x)=ky(\alpha x), y(x_0)=y_0,
\end{equation}
where $x_0\neq 0$, we have the following  results.

(i). If $x_0$ is not the zero of $E_{\alpha}(kx)$, that is, $E_{\alpha}(kx_0)\neq 0$, then there exists the unique solution
\begin{equation}
y(x)=\frac{y_0E_{\alpha}(kx)}{E_{\alpha}(kx_0)}.
\end{equation}

(ii). If $y_0=0$ and $x_0$ is the zero of $E_{\alpha}(kx)$, that is, $E_{\alpha}(kx_0)=0$, then there exist an infinity of solutions all of which can be represented by
\begin{equation}
y(x)=cE_{\alpha}(kx),
\end{equation}
where $c$ is an arbitrary constant.

(iii). If $y_0\neq 0$ and $E_{\alpha}(kx_0)=0$, then there does not exist solution.

\textbf{Proof}. For case (i), it is easy to see that $y(x)=\frac{y_0E_{\alpha}(kx)}{E_{\alpha}(kx_0)}$ is a solution for the initial value problem. Then we only need to prove that it is the unique solution. Assume that $y_1$ and $y_2$ are two solutions satisfying $y_1(x_0)=y_2(x_0)=y_0$. Letting $z(x)=y_1(x)-y_2(x)$, then we have
\begin{equation*}
z'(x)=kz(\alpha x), z(x_0)=0.
\end{equation*}
If we can prove $z(0)=0$, then we can get $z(x)\equiv 0$ by the unique theorem 2.1 at the origin point. In fact, if $z(0)=z_0\neq 0$, by the unique theorem 2.1 and proposition 3.1,  we must get
\begin{equation*}
z(x)=z_0E_{\alpha}(kx),
\end{equation*}
and hence
\begin{equation*}
z(x_0)=z_0E_{\alpha}(kx_0)\neq 0,
\end{equation*}
which is contradictory to $z(x_0)=0$.

For the case (ii), the existence of the solution is obvious. Let $y(x) $ be any solution with $y(x_0)=0$,  then $y(x)$ has a value $y(0)$ at $x=0$, and hence we have a solution $y(x)=y(0)E_{\alpha}(kx)$. By taking $y(0)=c$, we get the conclusion. The case (iii) is obvious. The proof is completed.

\textbf{Theorem 8.2}. For existence and uniqueness of the solution to the initial value problem of non-homogenous equation
\begin{equation}
y'(x)=\lambda y(\alpha x)+qE_{\alpha}(\lambda\alpha x),y(x_0)=y_0,
\end{equation}
there are the following three cases:

(i). If $E_{\alpha}(\lambda x_0)\neq 0$, then there exists the unique solution
\begin{equation}
y(x)=\frac{y_0-qx_0E_{\alpha}(\lambda\alpha x_0)}{E_{\alpha}(\lambda x_0)}E_{\alpha}(\lambda x)+qxE_{\alpha}(\lambda\alpha x).
\end{equation}

(ii). If $E_{\alpha}(\lambda x_0)=0$ and $y_0=qx_0E_{\alpha}(\lambda\alpha x_0)$, then there exist an infinity of solutions all of which can be given by
\begin{equation}
y(x)=cE_{\alpha}(\lambda x)+qxE_{\alpha}(\lambda\alpha x),
\end{equation}
where $c$ is an arbitrary constant.

(iii). If $E_{\alpha}(\lambda x_0)=0$ and $y_0\neq qx_0E_{\alpha}(\lambda\alpha x_0)$, then there does not exist solution.

\textbf{Proof}. (i). We only need to prove the uniqueness. Assuming that $z(x)$ is another solution with $z(0)=c$, by theorem 4.2, we have
\begin{equation*}
z(x)=cE_{\alpha}(\lambda x)+qxE_{\alpha}(\lambda\alpha x).
\end{equation*}
Since we also require $z(x_0)=y_0$, so we get $c=\frac{y_0-qx_0E_{\alpha}(\lambda\alpha x_0)}{E_{\alpha}(\lambda x_0)}$, and then $y(x)=z(x)$.

(ii). Firstly, we can verify that (87) is just solution and satisfies these two conditions. Hence we only need to prove that any solution satisfying these conditions has the form of (87).  Assume that $y(x)$ is the solution satisfying these two conditions. Then, according to the value $y(0)$ of $y(x)$ at $x=0$, we can use the theorem 4.2 to give the solution
\begin{equation*}
y(x)=y(0)E_{\alpha}(\lambda x)+qxE_{\alpha}(\lambda\alpha x).
\end{equation*}
Letting $y(0)=c$ gives conclusion.

(iii). It is obvious. The proof is completed.

\textbf{Theorem 8.3}. For the linear equations system
\begin{equation}
y'_1(x)=\lambda y_1(\alpha x)+y_{2}(\alpha x),
\end{equation}
\begin{equation}
y'_2(x)=\lambda y_2(\alpha x),
\end{equation}
with initial values $y_1(x_0)$ and $y_2(x_0)$,
 there are the following five cases about  its solution.

(i). If $E_{\alpha}(\lambda x_0)\neq 0$, then there exists the unique solution
\begin{equation}
y_1(x)=\frac{y_1(x_0)-\frac{y_2(x_0)}{E_{\alpha}(\lambda x_0)}x_0E_{\alpha}(\lambda\alpha x_0)}{E_{\alpha}(\lambda x_0)}E_{\alpha}(\lambda x)
+\frac{y_2(x_0)xE_{\alpha}(\lambda\alpha x)}{E_{\alpha}(\lambda x_0)},
\end{equation}
\begin{equation}
y_2(x)=\frac{y_2(x_0)}{E_{\alpha}(\lambda x_0)}E_{\alpha}(\lambda x).
\end{equation}

(ii). If $E_{\alpha}(\lambda x_0)=0$, $y_2(x_0)=0$ and $E_{\alpha}(\lambda\alpha x_0)\neq 0$, then there exist an infinity of solutions all of which can be given by
\begin{equation}
y_1(x)=cE_{\alpha}(\lambda x)+\frac{y_1(x_0)xE_{\alpha}(\lambda\alpha x)}{x_0E_{\alpha}(\lambda\alpha x_0)},
\end{equation}
\begin{equation}
y_2(x)=\frac{y_1(x_0)}{x_0E_{\alpha}(\lambda\alpha x_0)}E_{\alpha}(\lambda x),
\end{equation}
where $c$ is an arbitrary constant.

(iii). If $E_{\alpha}(\lambda x_0)=0$, $y_2(x_0)=0$ and $E_{\alpha}(\lambda\alpha x_0)=0$ and $y_1(x_0)=0$, then there exist an infinity of solutions all of which can be given by
\begin{equation}
y_1(x)=c_1E_{\alpha}(\lambda x)+c_2xE_{\alpha}(\lambda\alpha x),
\end{equation}
\begin{equation}
y_2(x)=c_2E_{\alpha}(\lambda x),
\end{equation}
where $c_1$ and $c_2$ are two arbitrary constants.

(iv). If $E_{\alpha}(\lambda x_0)=0$, $y_2(x_0)=0$ and $E_{\alpha}(\lambda\alpha x_0)=0$ and $y_1(x_0)\neq0$, then there does not exist solution.

(v). If $E_{\alpha}(\lambda x_0)= 0$ and $y_2(x_0)\neq0$, then there does not exist solutions.

\textbf{Proof}. By theorem 8.1 and theorem 8.2, we can easily prove it. The proof is completed.

 Now we consider the second order pantograph equation
\begin{equation}
y''(t)+py'(\alpha t)+qy(\alpha^2t)=0,
\end{equation}
with initial conditions
\begin{equation}
y(t_0)=A, y'(\frac{t_0}{\alpha})=B,
\end{equation}
where $t_0$ is a general point. Letting $x_1(t)=y(t)$ and $x_2(t)=y'(\frac{t}{\alpha})$, we transform the second order equation into the system of the first order equations  as follows
\begin{equation}
x_1'(t)=x_2(\alpha t),
\end{equation}
\begin{equation}
x_2'(t)=-\frac{q}{\alpha}x_1(\alpha t)-\frac{p}{\alpha}x_2(\alpha t),
\end{equation}
with initial conditions
\begin{equation}
x_1(t_0)=A, x_2(t_0)=B.
\end{equation}
For the coefficients matrix

\begin{equation}
K=\left(
  \begin{array}{ccc}
    0 & 1 \\
    -\frac{q}{\alpha} &  -\frac{p}{\alpha}\\
  \end{array}
\right),
\end{equation}
 when $p^2\neq 4q\alpha$, there exist two distinct eigenvalues $\lambda_1\neq\lambda_2$, then there is an invertible matrix $P$ such that
\begin{equation}
P^{-1}KP=\left(
  \begin{array}{ccc}
    \lambda_1 & 0 \\
   0 & \lambda_2\\
  \end{array}
\right).
\end{equation}
When $p^2= 4q\alpha$, we have  $\lambda_1=\lambda_2=\lambda=-\frac{p}{2\alpha}$, and the corresponding elementary factor is $(\lambda+\frac{p}{2\alpha})^2$, and hence there is an invertible matrix $P$ such that
\begin{equation}
P^{-1}KP=\left(
  \begin{array}{ccc}
    \lambda & 1 \\
   0 & \lambda\\
  \end{array}
\right).
\end{equation}
Denote

\begin{equation}
P=\left(
  \begin{array}{ccc}
    p_{11} & p_{12}\\
   p_{21} &p_{22}\\
  \end{array}
\right),
\end{equation}
\begin{equation}
P^{-1}=\left(
  \begin{array}{ccc}
    q_{11} & q_{12}\\
   q_{21} &q_{22}\\
  \end{array}
\right),
\end{equation}
and
\begin{equation}
x_1(t)=p_{11}z_1(t)+p_{12}z_2(t),
\end{equation}
\begin{equation}
x_2(t)=p_{21}z_1(t)+p_{22}z_2(t).
\end{equation}

According to the theorems 8.1-8.3, we can prove the following result.

\textbf{Theorem 8.4.} For the system of equations  (98) and (99) with condition (100), there are the following two cases to be discussed:

\emph{\textbf{Case 1}}. When $p^2\neq 4q\alpha$, then the matrix $K$ has two distinct eigenvalues $\lambda_1$ and $\lambda_2$, and hence we have
\begin{equation}
z'_1(t)=\lambda_1z_1(\alpha t),
\end{equation}
\begin{equation}
z'_2(t)=\lambda_2z_2(\alpha t),
\end{equation}
with initial conditions
\begin{equation}
z_1(t_0)=q_{11}x_1(t_0)+q_{12}x_2(t_0)=q_{11}A+q_{12}B,
\end{equation}
\begin{equation}
z_2(t_0)=q_{21}x_1(t_0)+q_{22}x_2(t_0)=q_{21}A+q_{22}B.
\end{equation}
Then there are the following six cases:

(i). If $E_{\alpha}(\lambda_1 t_0)\neq 0$ and $E_{\alpha}(\lambda_2 t_0)\neq 0$, then there exists the unique solution
\begin{equation}
y(t)=\frac{p_{11}(q_{11}A+q_{12}B)}{E_{\alpha}(\lambda_1 t_0)}E_{\alpha}(\lambda_1 t)
+\frac{p_{12}(q_{21}A+q_{22}B)}{E_{\alpha}(\lambda_2 t_0)}E_{\alpha}(\lambda_2 t),
\end{equation}
or equivalently
\begin{equation}
y(t)=\frac{A\lambda_2-B}{\lambda_2-\lambda_1}\frac{E_{\alpha}(\lambda_1 t)}{E_{\alpha}(\lambda_1 t_0)}
+\frac{A\lambda_1-B}{\lambda_1-\lambda_2}\frac{E_{\alpha}(\lambda_2 t)}{E_{\alpha}(\lambda_2 t_0)}.
\end{equation}

(ii). If $E_{\alpha}(\lambda_1 t_0)\neq 0$, $E_{\alpha}(\lambda_2 t_0)= 0$ and $z_2(t_0)=0$, then there exist an infinity of solutions all of which can be given by
\begin{equation}
y(t)=p_{11}\frac{q_{11}A+q_{12}B}{E_{\alpha}(\lambda_1 t_0)}E_{\alpha}(\lambda_1 t)
+p_{12}cE_{\alpha}(\lambda_2 t),
\end{equation}
or equivalently
\begin{equation}
y(t)=\frac{AE_{\alpha}(\lambda_1 t)}{E_{\alpha}(\lambda_1 t_0)}
+cE_{\alpha}(\lambda_2 t),
\end{equation}
where $c$ is an arbitrary constant.

(iii). If $E_{\alpha}(\lambda_2 t_0)\neq 0$, $E_{\alpha}(\lambda_1 t_0)= 0$ and $z_1(t_0)=0$, then there exist an infinity of solutions  all of which can be given by
\begin{equation}
y(t)=p_{12}\frac{q_{21}A+q_{22}B}{E_{\alpha}(\lambda_2 t_0)}E_{\alpha}(\lambda_2 t)
+p_{11}cE_{\alpha}(\lambda_1 t),
\end{equation}
or equivalently
\begin{equation}
y(t)=\frac{AE_{\alpha}(\lambda_2 t)}{E_{\alpha}(\lambda_2 t_0)}
+cE_{\alpha}(\lambda_1 t),
\end{equation}
where $c$ is an arbitrary constant.

(iv). If $E_{\alpha}(\lambda_1 t_0)= 0$ and $z_1(t_0)=0$, $E_{\alpha}(\lambda_2 t_0)= 0$ and $z_2(t_0)=0$, then $A=B=0$, and there exist an infinity of solutions  all of which can be given by
\begin{equation}
y(t)=c_1E_{\alpha}(\lambda_1 t)+c_2E_{\alpha}(\lambda_2 t),
\end{equation}
where $c_1$ and $c_2$ are two arbitrary constants.

(v). If $E_{\alpha}(\lambda_1 t_0)= 0$ and $z_1(t_0)\neq 0$, then there does not exist solution.

(vi). If  $E_{\alpha}(\lambda_2 t_0)= 0$ and $z_2(t_0)\neq0$, then there does not exist solution.

\emph{\textbf{Case 2}}. When $p^2=4q\alpha$, we have  $\lambda_1=\lambda_2=\lambda=-\frac{p}{2\alpha}$,  then there are the following five cases:

(i). If $E_{\alpha}(\lambda t_0)\neq 0$, then there exists the unique solution
\begin{equation*}
y(t)=p_{11}\{\frac{z_1(t_0)E_{\alpha}(\lambda t_0)-z_2(t_0)t_0E_{\alpha}(\lambda\alpha t_0)}{E^2_{\alpha}(\lambda t_0)}E_{\alpha}(\lambda t)
+\frac{z_2(t_0)tE_{\alpha}(\lambda\alpha t)}{E_{\alpha}(\lambda t_0)}\}+p_{12}\frac{z_2(t_0)E_{\alpha}(\lambda t)}{E_{\alpha}(\lambda t_0)},
\end{equation*}
or equivalently,
\begin{equation}
y(t)=\frac{AE_{\alpha}(\lambda t_0)-(B-A\lambda)t_0E_{\alpha}(\lambda\alpha t_0)}{E^2_{\alpha}(\lambda t_0)}E_{\alpha}(\lambda t)
+\frac{B-A\lambda}{E_{\alpha}(\lambda t_0)}tE_{\alpha}(\lambda\alpha t).
\end{equation}

(ii). If $E_{\alpha}(\lambda t_0)=0$, $E_{\alpha}(\lambda\alpha t_0)\neq 0$ and $z_2(t_0)=0$, then there exist an infinity of solutions all of which can be given by
\begin{equation*}
y(t)=p_{11}\{cE_{\alpha}(\lambda t)+\frac{z_1(t_0)tE_{\alpha}(\lambda\alpha t)}{t_0E_{\alpha}(\lambda\alpha t_0)}\}+
p_{12}\frac{z_1(t_0)E_{\alpha}(\lambda t)}{t_0E_{\alpha}(\lambda\alpha t_0)},
\end{equation*}
or equivalently,
\begin{equation}
y(t)=cE_{\alpha}(\lambda t)
+\frac{A}{t_0E_{\alpha}(\lambda\alpha t_0)}tE_{\alpha}(\lambda\alpha t),
\end{equation}
where $c$ is an arbitrary constant.

(iii). If $E_{\alpha}(\lambda t_0)= 0$, $E_{\alpha}(\lambda\alpha t_0)= 0$ and $z_1(t_0)=z_2(t_0)=0$, then $A=B=0$, and there exist an infinity of solutions  all of which can be given by
\begin{equation}
y(t)=c_1E_{\alpha}(\lambda t)
+c_2tE_{\alpha}(\lambda\alpha t),
\end{equation}
where $c_1$ and $c_2$ are two arbitrary constants.

(iv). If $E_{\alpha}(\lambda t_0)= 0$, $E_{\alpha}(\lambda\alpha t_0)= 0$ and $z_2(t_0)=0$ and $z_1(t_0)\neq0$, then there does not exist solutions.

(v). If $E_{\alpha}(\lambda t_0)= 0$ and $z_2(x_0)\neq0$, then there does not exist solutions.
.
\textbf{Proof}. We only consider (ii) in the case 1 and give a detailed proof of the equivalence of (110) and (111). Other cases can be proved similarly. Then, we only need to prove $p_{11}(q_{11}A+q_{12}B)=A$. In fact, since $\lambda_1$ is the first eigenvalue of the matrix $K$, we have $p_{21}=\lambda_1p_{11}$. Further, by $P^{-1}P=E$, it follows that $p_{11}q_{12}+p_{12}q_{22}=0$ and $ p_{11}q_{11}+p_{21}q_{12}=1$. Therefore, from $z_2(t_0)=0$, that is, $q_{21}A+q_{22}B=0$. we get
\begin{equation*}
\frac{B}{A}=-\frac{q_{21}}{q_{22}}=\frac{p_{21}}{p_{11}}=\lambda_1.
\end{equation*}
And hence, we have
\begin{equation*}
p_{11}(q_{11}A+q_{12}B)=A(p_{11}q_{11}+\lambda_1p_{11}q_{12})=A(p_{11}q_{11}+p_{21}q_{12})=A.
\end{equation*}
The proof is completed.

These above theorems are the basic results for the existence and uniqueness of the kind of pantograph equations. It is easy to generalize the theorems 8.3 and 8.4 to $n$ dimensional case, but the result is complicated, so we do not write it. By combining these results with the theorem 6.2 and theorem 7.3, we can easily give the corresponding existence and uniqueness  results for the initial value problems at a general point of linear pantograph equations system and high order linear pantograph equations. For simplicity, we also do not write them.  However, we must point out that at a general initial point, the solutions of the initial value problem are complicated, and includes no solution, unique solution and an infinity of solutions. All of these are rooted in the MFBH theorem. These results show the essential differences between the pantograph equation and usual ordinary differential equations.

\textbf{Remark 8.2}. We must have noticed that for the system of the second order linear pantograph equations (98)and (99), its initial conditions are very special, that is, conditions (100) of $y$ and $y'$ are taken respectively at two points $t_0$ and
$\frac{t_0}{\alpha}$ but not the only one point $t_0$. This is because we can use the equivalent linear pantograph equations system to deal with it. If the initial point is taken at only point $t_0$, we can use another method to give the following result.

\textbf{Theorem 8.5}. Consider the second order linear pantograph equation
\begin{equation}
y''(t)+py'(\alpha t)+qy(\alpha^2t)=0,
\end{equation}
with initial conditions
\begin{equation}
y(t_0)=A, y'(t_0)=B,
\end{equation}
where $t_0$ is a general point. Let $\lambda_1$ and $\lambda_2$ be two roots of its characteristic equation
\begin{equation}
\alpha\lambda^2+p\lambda+q=0.
\end{equation}
Then we have the following two cases:

\emph{\textbf{Case 1}}. $p^2\neq4\alpha q$, that is $\lambda_1\neq\lambda_2$. There are three cases for the solutions.

(i). If $\left|
\begin{array}{ccc}
E_{\alpha}(\lambda_1 t_0)& E_{\alpha}(\lambda_2 t_0)  \\
\lambda_1E_{\alpha}(\alpha\lambda_1 t_0) & \lambda_2E_{\alpha}(\alpha\lambda_2 t_0)
\end{array}
\right|\neq0$, then there is the unique solution
\begin{equation}
y(t)=c_1E_{\alpha}(\lambda_1 t)+c_2E_{\alpha}(\lambda_2 t),
\end{equation}
where $c_1$ and $c_2$ can be determined by the following equations system
\begin{equation}
c_1E_{\alpha}(\lambda_1 t_0)+c_2E_{\alpha}(\lambda_2 t_0)=A,
\end{equation}
\begin{equation}
c_1\lambda_1E_{\alpha}(\alpha\lambda_1 t_0)+c_2\lambda_2E_{\alpha}(\alpha\lambda_2 t_0)=B.
\end{equation}

(ii). If $\left|
\begin{array}{ccc}
E_{\alpha}(\lambda_1 t_0)& E_{\alpha}(\lambda_2 t_0)  \\
\lambda_1E_{\alpha}(\alpha\lambda_1 t_0) & \lambda_2E_{\alpha}(\alpha\lambda_2 t_0)
\end{array}
\right|=0$, and the rank of the matrix
$\left(
  \begin{array}{ccc}
    E_{\alpha}(\lambda_1 t_0)& E_{\alpha}(\lambda_2 t_0)& A  \\
\lambda_1E_{\alpha}(\alpha\lambda_1 t_0) & \lambda_2E_{\alpha}(\alpha\lambda_2 t_0)& B
  \end{array}
\right)$
is 1, then there is an infinity of solutions  all of which can be given by
\begin{equation}
y(t)=c_1E_{\alpha}(\lambda_1 t)+c_2E_{\alpha}(\lambda_2 t),
\end{equation}
where $c_1$ and $c_2$ satisfy (126) and (127) which has an infinity of solutions.

(iii). If $\left|
\begin{array}{ccc}
E_{\alpha}(\lambda_1 t_0)& E_{\alpha}(\lambda_2 t_0)  \\
\lambda_1E_{\alpha}(\alpha\lambda_1 t_0) & \lambda_2E_{\alpha}(\alpha\lambda_2 t_0)
\end{array}
\right|=0$, and the rank of the matrix
$\left(
  \begin{array}{ccc}
    E_{\alpha}(\lambda_1 t_0)& E_{\alpha}(\lambda_2 t_0)& A  \\
\lambda_1E_{\alpha}(\alpha\lambda_1 t_0) & \lambda_2E_{\alpha}(\alpha\lambda_2 t_0)& B
  \end{array}
\right)$
is 2, then there does not exist solution.

\emph{\textbf{Case 2}}.  $p^2=4\alpha q$, that is $\lambda_1=\lambda_2=\lambda$. There are three cases for the solutions.

(i). If $\left|
\begin{array}{ccc}
E_{\alpha}(\lambda t_0)& t_0E_{\alpha}(\lambda\alpha t_0)  \\
\lambda E_{\alpha}(\alpha\lambda t_0) & E_{\alpha}(\alpha\lambda t_0)+\lambda\alpha t_0E_{\alpha}(\alpha^2\lambda t_0)
\end{array}
\right|\neq0$, then there is the unique solution
\begin{equation}
y(t)=c_1E_{\alpha}(\lambda t)+c_2tE_{\alpha}(\lambda \alpha t),
\end{equation}
where $c_1$ and $c_2$ can be determined by the following equations system
\begin{equation}
c_1E_{\alpha}(\lambda t_0)+c_2t_0E_{\alpha}(\lambda\alpha t_0)=A,
\end{equation}
\begin{equation}
c_1\lambda E_{\alpha}(\alpha\lambda t_0)+c_2(E_{\alpha}(\alpha\lambda t_0)+\lambda\alpha t_0 E_{\alpha}(\alpha^2\lambda t_0))=B.
\end{equation}

(ii). If $\left|
\begin{array}{ccc}
E_{\alpha}(\lambda t_0)& t_0E_{\alpha}(\lambda\alpha t_0)  \\
\lambda E_{\alpha}(\alpha\lambda t_0) & E_{\alpha}(\alpha\lambda t_0)+\lambda\alpha t_0E_{\alpha}(\alpha^2\lambda t_0)
\end{array}
\right|=0$, and the rank of the matrix $\left(
  \begin{array}{ccc}
   E_{\alpha}(\lambda t_0)& t_0E_{\alpha}(\lambda\alpha t_0)& A  \\
\lambda E_{\alpha}(\alpha\lambda t_0) & E_{\alpha}(\alpha\lambda t_0)+\lambda\alpha t_0E_{\alpha}(\alpha^2\lambda t_0)& B
  \end{array}
\right)$ is 1, then there is an infinity of solutions  all of which can be given by
\begin{equation}
y(t)=c_1E_{\alpha}(\lambda_1 t)+c_2E_{\alpha}(\lambda_2 t),
\end{equation}
where $c_1$ and $c_2$ satisfy (130) and (131) which has an infinity of solutions.

(iii). If $\left|
\begin{array}{ccc}
E_{\alpha}(\lambda t_0)& t_0E_{\alpha}(\lambda\alpha t_0)  \\
\lambda E_{\alpha}(\alpha\lambda t_0) & E_{\alpha}(\alpha\lambda t_0)+\lambda\alpha t_0E_{\alpha}(\alpha^2\lambda t_0)
\end{array}
\right|=0$, and the rank of the matrix $\left(
  \begin{array}{ccc}
   E_{\alpha}(\lambda t_0)& t_0E_{\alpha}(\lambda\alpha t_0)& A  \\
\lambda E_{\alpha}(\alpha\lambda t_0) & E_{\alpha}(\alpha\lambda t_0)+\lambda\alpha t_0E_{\alpha}(\alpha^2\lambda t_0)& B
  \end{array}
\right)$ is 2, then there does not exist solution.

\textbf{Proof}. We only prove the case 1, the case 2 can be proven similarly. In fact, since the solution and its derivative can take values at the point $t=0$, we can get the general solution as
\begin{equation*}
y(t)=c_1E_{\alpha}(\lambda_1 t)+c_2E_{\alpha}(\lambda_2 t).
\end{equation*}
Therefore, we have
\begin{equation*}
A=y(t_0)=c_1E_{\alpha}(\lambda_1 t_0)+c_2E_{\alpha}(\lambda_2 t_0),
\end{equation*}
\begin{equation*}
B=y'(t_0)=c_1\lambda_1E_{\alpha}(\lambda_1\alpha t_0)+c_2\lambda_2E_{\alpha}(\lambda_2\alpha t_0).
\end{equation*}
According to the theory of linear algebraic equations system, we get the corresponding conclusions (i-iii) of case 1. The proof is completed.

\textbf{Remark 8.3}. We have seen that for general linear homogenous pantograph equations, we can write their general solutions. When we consider the initial value problem at a general point,  we will need to determine these constants in general solutions. However, since exponent-like function $E_{\alpha}(x)$ has an infinity of real zeroes, we can't find these coefficients or can find an infinite number of solutions in some cases. If we further consider the existence, uniqueness and non-uniqueness of the solution of the initial value problem at a general point for the non-homogenous linear pantograph equations, we also obtain the similar results by the similar consideration.

 In the theorem 8.2. we have considered the first  order non-homogenous linear pantograph equation with the resonance term. Next, we give the result in the non-resonance case. Other cases such as the second order linear non-homogenous pantograph equation can be dealt with by the similar method, and the results can be derived from the corresponding results on the homogenous equations such as the theorems 8.5 and 7.5.

\textbf{Theorem 8.6}. For the non-homogenous pantograph equation
\begin{equation*}
y'(x)+\beta y(\alpha x)=AE_{\alpha}(rx),y(x_0)=y_0,
\end{equation*}
where $x_0\neq0$ and $r\neq\alpha\beta$. There are the following three cases:

(i). If $E_{\alpha}(\beta x_0)\neq0$, then there exists the unique solution
\begin{equation}
y(x)=\frac{y_0-\frac{A\alpha}{r-\alpha\beta}E_{\alpha}(\frac{r}{\alpha} x_0)}{E_{\alpha}(\beta x_0)}E_{\alpha}(\beta x)+\frac{A\alpha}{r-\alpha\beta}E_{\alpha}(\frac{r}{\alpha} x).
\end{equation}

(ii). If $E_{\alpha}(\beta x_0)=0$ and $y_0=\frac{A\alpha}{r-\alpha\beta}E_{\alpha}(\frac{r}{\alpha} x_0)$, then there exist an infinity of solutions all of which can be given by
\begin{equation}
y(x)=cE_{\alpha}(\beta x)+\frac{A\alpha}{r-\alpha\beta}E_{\alpha}(\frac{r}{\alpha} x),
\end{equation}
where $c$ is an arbitrary constant.

(iii). If $E_{\alpha}(\beta x_0)=0$ and $y_0\neq\frac{A\alpha}{r-\alpha\beta}E_{\alpha}(\frac{r}{\alpha} x_0)$, then there does not exist the solution.

\textbf{Proof}. According to the general solution (67), we can easily prove it. The proof is completed.

Finally, we give a generalization of the MFBH theorem to the first order linear pantograph equation with variable coefficient. By the similar method of proving MFBH theorem by Hahn[8], we give the proofs of the following two results.

\textbf{Theorem 8.7}. For the equation (here $0<\alpha<1$)
\begin{equation}
y'(t)=-k(t)y(\alpha t), y(0)\neq0,
\end{equation}
if $k'(t)\leq0$ and $k(t)>k_0$ where $k_0>0$ is a constant, then the solution has an infinity of positive zeros.

\textbf{Proof}. We will prove that it is impossible that $y(t)>0$ for all $t>0$. Suppose that  $y(t)>0$ for all $t>0$, then $y(t)>0$ for $t>\alpha^2t_0$ where $t_0>0$ is an arbitrary constant. Therefore, for $t>t_0$, we have

\begin{equation*}
y(t)>0, y'(t)<0,
\end{equation*}
\begin{equation*}
y''(t)=-k'(t)y(\alpha t)+\alpha k(t)k(\alpha t)y(\alpha^2 t)>0.
\end{equation*}
Take a sequence of points $t_n=\frac{t_0}{\alpha^n}$ for $n=0,1,\cdots.$ On the interval $[t_n,t_{n+1}]$, the graph of the function $y(t)$ is convex and decreasing, so the arc $A_nA_{n+1}$ lies above the tangent line $L$ at point $A_{n+1}=(t_{n+1},y(t_{n+1}))$ and below the line through the point $A_n=(t_n,y(t_n))$ parallel to $L$. And then, we have
\begin{equation*}
y(t_{n+1})<y(t_n)(1-k(t_n)(t_{n+1}-t_n))=y(t_n)(1-\frac{k(t_n)}{\alpha^n}(t_{1}-t_0)).
\end{equation*}
Since $k(t)>k_0>0$ and $0<\alpha<1$, we have $\frac{k(t_n)}{\alpha^n}\rightarrow+\infty$ as $n\rightarrow+\infty$, and hence $1-\frac{k(t_n)}{\alpha^n}(t_{1}-t_0)$ will become negative for sufficiently large $n$. This is a contradiction with the assumption.  By the similar reason, $y(t)$ can not be negative for sufficiently large $t>0$. The proof is completed.

Similarly, we have the following result.

\textbf{Theorem 8.8}. For the equation (here $0<\alpha<1$)
\begin{equation}
y'(t)=k(t)y(\alpha t), y(0)\neq0,
\end{equation}
if $k'(t)\geq0$ and $k(t)>k_0$ where $k_0>0$ is a constant, then every solution has an infinity of negative zeros.

\section{Boundary value problem of the second order linear pantograph equation and its applications}

\subsection{Boundary value problem of the second order linear pantograph equation}
 Consider the boundary value problem
\begin{equation}
y''(x)=\lambda y(\alpha^2 x), (0<\alpha<1),
\end{equation}
with boundary condition
\begin{equation}
y(0)=y(1)=0.
\end{equation}
If $\alpha=1$, then we have $\lambda<0$ for nontrivial solution. However, if $0<\alpha<1$, we need some efforts to proof this result.

\textbf{Lemma 9.1}. The function $h(x)=E_{\alpha}(x)-E_{\alpha}(-x)$ has the unique real zero $x=0$.

\textbf{Proof}. We can give a direct proof by the Taylor expansion of $E_{\alpha}(x)$. Here we give another proof. First, since $h(-x)=E_{\alpha}(-x)-E_{\alpha}(x)=-h(x)$, we know that if $h(x_0)=0$, then $h(-x_0)=0$. This means that if there is a positive zero, then there is a corresponding negative zero. Therefore, we only need to prove that $h(x)$ has no positive zero. In fact, $h(x)$ satisfies the following equation
\begin{equation}
h''(x)=\alpha h(\alpha^2 x),
\end{equation}
with initial conditions $h(0)=0$ and $h'(0)=2$. By the theorem 2.1, there exists the unique analytic solution. Thus, we can assume that the expansion of $h(x)$ is
\begin{equation}
h(x)=\sum_{n=0}^{+\infty}b_nx^n,
\end{equation}
where $b_n's$ are the undetermined coefficients. Substituting it into equation gives
\begin{equation}
b_{n+2}=\frac{b_n\alpha^{2n+1}}{(n+2)(n+1)}, n=0,1,2,\cdots.
\end{equation}
Since $b_0=0$ and $b_1=2$, so $b_{2n}=0$ and $b_{2n-1}>0$ for $n>0$,  and then $h(x)$ has no positive zero. The proof is completed.

Similarly, we have the following lemma.

\textbf{Lemma 9.2}. The function $h(x)=E_{\alpha}(x)+E_{\alpha}(-x)$ has no any real zero.

By the above lemma 9.1, we have the following theorem.

\textbf{Theorem 9.1}. For the boundary value problem (137) and (138), if there is nontrivial solution, we must have  $\lambda<0$.

\textbf{Proof}. Assume that $\lambda>0$. Then the general solution of equation (137) is give by
\begin{equation}
y(x)=c_1E_{\alpha}(\sqrt {\frac{\lambda}{\alpha}} x)+c_2E_{\alpha}(-\sqrt {\frac{\lambda}{\alpha}} x).
\end{equation}
From the boundary conditions (141), we have
\begin{equation}
c_1+c_2=0,
\end{equation}
\begin{equation}
c_1E_{\alpha}(\sqrt {\frac{\lambda}{\alpha}})+c_2E_{\alpha}(-\sqrt {\frac{\lambda}{\alpha}})=0.
\end{equation}
By the lemma 9.1, the coefficients determinant is not zero, that is,
\begin{equation}
\left|
\begin{array}{ccc}
1 & 1 \\
 E_{\alpha}(\sqrt {\frac{\lambda}{\alpha}}) & E_{\alpha}(-\sqrt {\frac{\lambda}{\alpha}})
\end{array}
\right|\neq0,
\end{equation}
so $c_1=c_2=0$, and then the solution $y(x)\equiv 0$ is trivial. The proof is completed.

\textbf{Theorem 9.2}. For the boundary value problem (137) and (138), there exist an infinity of negative eigenvalues $\lambda_n (n=1,2,\cdots)$ which satisfy
\begin{equation}
\lambda_n=-\alpha\rho_n^2,
\end{equation}
and every corresponding eigenfunction is given by
\begin{equation}
y_n(x)=S_{\alpha}(\rho_n x),
\end{equation}
where $\rho_n$ is the $n$-th positive zero of $S_{\alpha}(x)$.

\textbf{Proof}. By theorem 9.1, we must have $\lambda<0$, then the general solution of the second order pantograph equation (137) is given by
\begin{equation}
y(x)=c_1C_{\alpha}(\sqrt {-\frac{\lambda}{\alpha}} x)+c_2S_{\alpha}(\sqrt {-\frac{\lambda}{\alpha}} x).
\end{equation}
From $y(0)=y(1)=0$, we have $c_1=0$ and
\begin{equation}
c_2S_{\alpha}(\sqrt {-\frac{\lambda}{\alpha}})=0.
\end{equation}
Since $c_2\neq 0$, $\lambda$ must satisfies
\begin{equation}
S_{\alpha}(\sqrt {-\frac{\lambda}{\alpha}})=0.
\end{equation}
Denote all zeroes of $S_{\alpha}(x)$ as $\rho_0=0, \pm\rho_1,\pm\rho_2,\cdots,$. So we have
an infinity of $\lambda's$ as follows,
\begin{equation}
\lambda_n=-\alpha\rho^2_n, n=1,2,\cdots.
\end{equation}
Therefore, for each eigenvalue $\lambda_n$, an eigenfunction is given by
\begin{equation*}
y_n(x)=S_{\alpha}(\rho_n x).
\end{equation*}
The proof is completed.

\textbf{Theorem 9.3}. For the symmetric boundary condition
\begin{equation}
y(-l)=y(l)=0,
\end{equation}
the eigenvalues of equation (137) must be negative and are given by
\begin{equation}
\lambda_{2n}=-\alpha\frac{\rho^2_n}{l^2}, \lambda_{2n-1}=-\alpha\frac{\eta^2_n}{l^2}, n=1,2,\cdots,
\end{equation}
and corresponding basic eigenfunctions are
\begin{equation}
y_{2n}(x)=S_{\alpha}(\frac{\rho_n}{l}x),
\end{equation}
\begin{equation}
y_{2n+1}(x)=C_{\alpha}(\frac{\eta_n}{l}x),
\end{equation}
where $\rho_n$ and $\eta_n$ are respectively the $n$-th zeros of $S_{\alpha}(x)$ and $C_{\alpha}(x)$.

\textbf{Proof}. By lemmas 9.1 and 9.2, if there exists nontrivial solution, we must have $\lambda<0$. In fact, if we suppose $\lambda>0$,  from the boundary conditions (152) and the general solution (142), we have
\begin{equation}
c_1E_{\alpha}(\sqrt {\frac{\lambda}{\alpha}}l)+c_2E_{\alpha}(-\sqrt {\frac{\lambda}{\alpha}}l)=0
\end{equation}
\begin{equation}
c_1E_{\alpha}(-\sqrt {\frac{\lambda}{\alpha}}l)+c_2E_{\alpha}(\sqrt {\frac{\lambda}{\alpha}}l)=0.
\end{equation}
The coefficients determinant is
\begin{equation}
\left|
\begin{array}{ccc}
E_{\alpha}(\sqrt {\frac{\lambda}{\alpha}}l) & E_{\alpha}(-\sqrt {\frac{\lambda}{\alpha}}l) \\
 E_{\alpha}(-\sqrt {\frac{\lambda}{\alpha}}l) & E_{\alpha}(\sqrt {\frac{\lambda}{\alpha}}l)
\end{array}
\right|=E^2_{\alpha}(\sqrt {\frac{\lambda}{\alpha}}l)-E^2_{\alpha}(-\sqrt {\frac{\lambda}{\alpha}}l).
\end{equation}
By the lemmas 9.1 and 9.2, the determinant is not equal to zero, so $c_1=c_2=0$, and hence the solution $y(x)\equiv 0$ is trivial. Therefore, we must have $\lambda<0$, and the general solution is (148). By boundary conditions (152), we have
\begin{equation}
c_1C_{\alpha}(\sqrt {-\frac{\lambda}{\alpha}} l)+c_2S_{\alpha}(\sqrt {-\frac{\lambda}{\alpha}} l)=0
\end{equation}
\begin{equation}
c_1C_{\alpha}(\sqrt {-\frac{\lambda}{\alpha}} l)-c_2S_{\alpha}(\sqrt {-\frac{\lambda}{\alpha}} l)=0.
\end{equation}
Since $c_1$ and $c_2$ can not be all zero, it follows that the determinant of the coefficients is zero, that is,
\begin{equation}
C_{\alpha}(\sqrt {-\frac{\lambda}{\alpha}} l)S_{\alpha}(\sqrt {-\frac{\lambda}{\alpha}} l)=0,
\end{equation}
whose solutions are
\begin{equation*}
\lambda_{2n}=-\alpha\frac{\rho^2_n}{l^2}, \lambda_{2n-1}=-\alpha\frac{\eta^2_n}{l^2}, n=1,2,\cdots.
\end{equation*}
The corresponding eigenfunctions can be taken as
\begin{equation*}
y_{2n}(x)=S_{\alpha}(\frac{\rho_n}{l}x),
\end{equation*}
\begin{equation*}
y_{2n-1}(x)=C_{\alpha}(\frac{\eta_n}{l}x),
\end{equation*}
for $n=1,2,\cdots$. The proof is completed.

\textbf{Remark 9.1}. We must notice that these eigenfunctions are not orthogonal each other.

\subsection{Formal solutions of heat-like pantograph equation and wave-like equation}

 Finally, we give an application of the results in section 9.1. These are only formal derivations. Consider the following  heat-like pantograph equation
\begin{equation}
u_t(\alpha^2 x,t)=u_{xx}(x,\beta t),
\end{equation}
\begin{equation}
u(0,t)=u(1,t)=0,
\end{equation}
\begin{equation}
u(x,0)=\phi(x).
\end{equation}
If $\alpha=\beta=1$, it will reduce to the usual heat equation.

We use the method of variables separation to solve the heat-like pantograph equation. Letting $u(x,t)=X(x)T(t)$ and
substituting it into the equation (162) yields two equations
\begin{equation}
T'(t)=\lambda T(\beta t),
\end{equation}
and
\begin{equation}
X''(x)=\lambda X(\alpha^2 x),
\end{equation}
with boundary condition
\begin{equation*}
X(0)=X(1)=0.
\end{equation*}

By the theorems 9.1 and 9.2, we have $\lambda_n=-\alpha\rho^2_n$ for $n=1,2,\cdots$ and
\begin{equation}
X_n(x)=S_{\alpha}(\rho_n x),
\end{equation}
and correspondingly
\begin{equation}
T_n(t)=E_{\beta}(-\alpha\rho^2_n t).
\end{equation}
Thus, formally, we have
\begin{equation}
u(x,t)=\sum_{n=1}^{+\infty}A_nE_{\beta}(-\alpha\rho^2_n t)S_{\alpha}(\rho_n x),
\end{equation}
which satisfies the heat-like pantograph  equation and the boundary condition. We now use the initial condition to
 solve the coefficients $A_n$. In fact, taking $t=0$ in the above formula gives
\begin{equation}
\phi(x)=\sum_{n=1}^{+\infty}A_nS_{\alpha}(\rho_n x),
\end{equation}
that is, $\phi(x)$ can be expanded as a Fourier's series based on sine-like functions. Unfortunately,
$\{S_{\alpha}(\rho_n x)\}_{n=1}^{\infty}$ is not a system of orthogonal bases, and hence we can not direct
compute the value of $A_n$. However, we can use the Gram-Schimdt's process to compute these coefficients.
Denote $f_n=S_{\alpha}(\rho_n x)$ for $n=1,2,\cdots$ and $\{e_n\}_{n=1}^{+\infty}$ the corresponding
orthogonal bases given by Gram-Schimdt's process. Then, for $n=2,3,\cdots$, we have
\begin{equation*}
e_1=f_1,
\end{equation*}
\begin{equation}
e_n=f_n-\frac{<f_n,e_1>}{<e_1,e_1>}e_1-\cdots-\frac{<f_n,e_{n-1}>}{<e_{n-1},e_{n-1}>}e_{n-1},
\end{equation}
where inner product $<f,g>=\int_0^1f(x)g(x)\mathrm{d}x$. Then $e_n$ is the linear combination of $f_1, \cdots, f_n$,
that is,
\begin{equation}
e_n=\sum_{m=1}^{n}H_{nm}f_m,
\end{equation}
where $H_{nm}$ can be explicitly  represented in terms of $\frac{<f_i,e_{j}>}{<e_{j},e_{j}>}$ for $i,j=1,\cdots,n.$
Therefore, we get
\begin{equation}
\phi(x)=\sum_{m=1}^{+\infty}A_mf_m=\sum_{n=1}^{+\infty}B_ne_n=\sum_{n=1}^{+\infty}B_n\sum_{m=1}^{n}H_{nm}f_m
=\sum_{m=1}^{+\infty}\sum_{n=m}^{+\infty}B_nH_{nm}f_m.
\end{equation}
It follows that
\begin{equation}
A_m=\sum_{n=m}^{+\infty}B_nH_{nm},
\end{equation}
and
\begin{equation}
B_n=\frac{<\phi(x),e_n>}{<e_n,e_n>}=\frac{\sum_{k=1}^{n}<\phi,f_k>}{\sum_{k=1}^{n}H^2_{nk}<f_k,f_k>},
\end{equation}
and hence, we have
\begin{equation}
A_m=\sum_{n=m}^{+\infty}\frac{H_{nm}\sum_{k=1}^{n}<\phi,f_k>}{\sum_{k=1}^{n}H^2_{nk}<f_k,f_k>}.
\end{equation}
Therefore, we give the formal solution of the heat-like pantograph equation
\begin{equation}
u(x,t)=\sum_{n=1}^{+\infty}\sum_{m=n}^{+\infty}\frac{H_{mn}\sum_{k=1}^{m}<\phi,f_k>}
{\sum_{k=1}^{m}H^2_{jk}<f_k,f_k>} E_{\beta}(-\alpha\rho^2_n
t)S_{\alpha}(\rho_n x).
\end{equation}

Similarly, we consider the following wave-like pantograph equation
\begin{equation}
u_{tt}(\alpha^2 x,t)=u_{xx}(x,\beta^2 t),
\end{equation}
\begin{equation}
u(0,t)=u(1,t)=0,
\end{equation}
\begin{equation}
u(x,0)=\phi(x),
\end{equation}
\begin{equation}
u_t(x,0)=\varphi(x),
\end{equation}
and give its formal solution
\begin{equation}
u(x,t)=\sum_{m=1}^{+\infty}\{\sum_{n=m}^{+\infty}\frac{H_{nm}\sum_{k=1}^{n}<\phi,f_k>}{\sum_{k=1}^{n}H^2_{nk}<f_k,f_k>}
C_{\beta}(-\rho_m t)
\end{equation}
\begin{equation}
-\frac{1}{\rho_m^2}\sum_{n=m}^{+\infty}\frac{H_{nm}\sum_{k=1}^{n}<\varphi,f_k>}{\sum_{k=1}^{n}H^2_{nk}<f_k,f_k>}
S_{\beta}(-\rho_m t)\}S_{\alpha}(\rho_m x).
\end{equation}

If we consider the wave-like pantograph equation on the infinite interval $(-\infty, +\infty)$,

\begin{equation}
u_{tt}(\alpha^2 x,t)=u_{xx}(x,\beta^2 t),
\end{equation}
\begin{equation}
u(-\infty,t)=u(+\infty,t)=0,
\end{equation}
\begin{equation}
u(x,0)=\phi(x),
\end{equation}
\begin{equation}
u_t(x,0)=\varphi(x),
\end{equation}
we will give the formal solution
\begin{equation*}
u(x,t)=\int_{-\infty}^{+\infty}S_{\beta}(\frac{y}{\sqrt\beta}t)
\{A_1(y)S_{\alpha}(\frac{y}{\sqrt\alpha}x)+B_1(y)C_{\alpha}(\frac{y}{\sqrt\alpha}x)\}\mathrm{d}y
\end{equation*}
\begin{equation}
+\int_{-\infty}^{+\infty}C_{\beta}(\frac{y}{\sqrt\beta}t)
\{A_2(y)S_{\alpha}(\frac{y}{\sqrt\alpha}x)+B_2(y)C_{\alpha}(\frac{y}{\sqrt\alpha}x)\}\mathrm{d}y,
\end{equation}
 where $A_k(y)$ and $B_k(y)$ ($k=1,2$) satisfy
\begin{equation}
\phi(x)=\int_{-\infty}^{+\infty}\{A_1(y)S_{\alpha}(\frac{y}{\sqrt\alpha}x)+B_1(y)C_{\alpha}(\frac{y}{\sqrt\alpha}x)\}\mathrm{d}y,
\end{equation}
and
\begin{equation}
\varphi(x)=\int_{-\infty}^{+\infty}\frac{y}{\sqrt\beta}\{A_2(y)S_{\alpha}(\frac{y}{\sqrt\alpha}x)+B_2(y)C_{\alpha}(\frac{y}{\sqrt\alpha}x)\}\mathrm{d}y.
\end{equation}

\textbf{Remark 9.2.}  We can formally introduce the Fourier-like transformation
 \begin{equation}
f(x)=\int_{-\infty}^{+\infty}F(y)E_{\alpha}(iyx)\mathrm{d}y,
\end{equation}
and the Laplace-like
transformation
\begin{equation}
L(f)(p)=\int_{0}^{+\infty}E_{\alpha}(-px)f(x)\mathrm{d}x,
\end{equation}
but the inverse transformations are unknown.

\textbf{Acknowledgments}: Thanks to Mr. Ballstadt for his pointing out the
papers[5,6]. I am also grateful to Prof. Iserles  for his kindly answers to my questions. Finally, thanks to Yue Kai for his helpful discussions.


\begin{thebibliography}{2}
\bibitem{w1}J K Hale. Functional differential equations. Springer, Berlin Heidelberg, 1971.
\bibitem{12}K Mahler. On a special functional equation. J. london
Math. Soc. 1940,1(2):115-123.
\bibitem{10}L Fox, D F Mayers, J R Ockendon and A B Tayler. On a
functional differantical equation.  IMA Journal of Applied Mathematics, 1971, 8(3): 271-307.
\bibitem{11}J R Ockendon and A B Tayler. The dynamics of a current
collection system for an electric locomotive. Proceedings of the Royal Society of London A: Mathematical, Physical and Engineering Sciences. The Royal Society, 1971, 322(1551): 447-468.
\bibitem{1}T Kato and J B McLeod. The functional-differential equation $y'(x)=ay(\lambda x)+by(x)$.
Bull. Amer. Math. Soc. 1971, 77:891-937.
\bibitem{2}J Carr and J Dyson. The functional differential equation $y'(x)=ay(\lambda x)+by(x)$. Proc.
Roy. Soc. Edinburgh Sect. A. 1974-75, 74: 165-174
\bibitem{3}J Carr and J Dyson. The matrix functional differential Equation $y'(x)=Ay(\lambda x)+By(x)$.
Proc. Roy. Soc. Edinburgh Sect. A.1975-76, 75: 5-22.
\bibitem{6}G R Morris, A Feldstein, and E W Bowen, The Phragmen¡ä ¡§Lindelof principle and a
class of functional differential equations. ¡®¡®Ordinary
Differential Equations,¡¯¡¯ pp. 513-540, Academic Press, New York,
1972.
\bibitem{8} A Iserles. On the generalized pantograph functional-differential equation. European J.
Appl. Math. 1992, 4(01): 1-38.
\bibitem{7}R Nussbaum. Existence and uniqueness theorems for some functional differential
equations of neutral type. Journal of Differential Equations.
1972,11(3):607-623.

\bibitem{4}A Iserles. On nonlinear delay differential equations. Trans. Amer. Math. Soc.
1994, 344(1): 441-477.
\bibitem{16}G Derfel, A Iserles. The pantograph equation in the complex plane.
 Journal of Mathematical Analysis and Applications, 1997, 213(1): 117-132.
 \bibitem{17}Yunkang Liu. Regular solutions of the Shabat equation. Journal of Differential Equations,
 1999, 154(1): 1-41.
\bibitem{5}Yunkang Liu. Asymptotic behaviour of functional-differential equations with proportional time
delays. European J. Appl. Math. 1996, 7(01):11-30.
\bibitem{24}J Cermak, P Kundrat, M Urbanek. Delay equations on time scales:
 Essentials and asymptotics of the solutions.
Journal of Difference Equations and Applications, 2008, 14(6): 567-580.
\bibitem{9} Yang Kuang, A Feldstein. Monotonic and oscillatory
solution of a linear neutral delay equation with infinite lag. SIAM
J. Math. Anal. 1990, 21(6): 1633-1641.
\bibitem{19}Yunkang Liu. Stability analysis of\ theta-methods for neutral functional-differential equations.
Numerische Mathematik, 1995, 70(4): 473-485.
\bibitem{20}A Iserles, Yunkang Liu. On pantograph integro-differential equations.
 University of Cambridge, Department of Applied Mathematics and Theoretical Physics, 1993.
\bibitem{23}J Mallet-Paret, R D Nussbaum. Analyticity and nonanalyticity of solutions of
 delay-differential equations.
 SIAM Journal on Mathematical Analysis, 2014, 46(4): 2468-2500.


\bibitem{14}Yunkang Liu. On Some Conjectures by Morriset al. about Zeros of an Entire Function.
Journal of Mathematical Analysis and Applications, 1998, 226(1):
1-5.
\bibitem{22}B Van Brunt, G C Wake. A Mellin transform solution to
 a second-order pantograph equation with linear dispersion arising in a cell growth model.
 European Journal of Applied Mathematics, 2011, 22(02): 151-168.
\bibitem{25}A Iserles, Yunkang Liu. Integro-differential equations and generalized hypergeometric functions.
Department of Applied Mathematics and Theoretical Physics, University of Cambridge, 1995.
\bibitem{g}A Feldstein, A Iserles, D Levin. Embedding of delay equations into an infinite-dimensional ODE system. Journal of Differential Equations, 1995, 117(1): 127-150.
\bibitem{w2}M Atiyah, G W Moore. A shifted view of fundamental physics. Surveys in Differential Geometry, 2010, 15.
arXiv:1009.3176v1.
\bibitem{w3}D X Kong, C Zhang. A new kind of functional differential equations. arXiv preprint arXiv:1402.3084, 2014.

\bibitem{2}V Spiridonov. Universal superpositions of coherent states and self-similar potentials. Physical Review A, 1995, 52(3): 1909.
\bibitem{3}S Skorik, V Spiridonov. Self-similar potentials and the q-oscillator algebra at roots of unity. Letters in Mathematical Physics, 1993, 28(1): 59-74.
\bibitem{1}L V Bogachev, G Derfel and S A Molchanov. On bounded continuous solutions of the archetypal equation with rescaling. Proc. R. Soc. A. The Royal Society, 2015, 471(2180): 20150351.
 \bibitem{el}Hill R. Mathematical theory of plasticity. Claredon
Press. 1950.
\bibitem{gu}V V Golubev. Lectures on the analytic theory of differential equations. Gostekhizdat, Moscow, 1950.
\bibitem{tm}E C Titchmarsh. The theory of functions.  London: Oxford University Press, 1952.
\bibitem{ps}G P\'{o}lya, J Schur.  \"{U}ber zwei Arten von Faktorenfolgen in der Theorie der algebraischen Gleichungen. Journal f\"{u}r die reine und angewandte Mathematik. 1914, 144: 89-113.
\bibitem{aaa}E Hille. Analytic function theory, Volume II. The second Edition. Providence Rhode Island: American Mathematical Soc., 2005.
\bibitem{Cl}Cheng-shi Liu. The asymptotic formulas of zeros of solutions for $y''(x)=-\alpha y(\alpha^2 x).$ In prepearation. 
\end{thebibliography}
\end{document}